\documentclass[11pt]{article}
\usepackage[utf8]{inputenc}
\usepackage[T1]{fontenc}
\usepackage{amsmath,amssymb,amsthm}
\usepackage[margin=1in]{geometry}
\usepackage{tikz}
\usepackage{hyperref}
\usepackage{xcolor}
\newtheorem{theorem}{Theorem}[section]
\newtheorem{lemma}[theorem]{Lemma}
\newtheorem{proposition}[theorem]{Proposition}
\newtheorem{corollary}[theorem]{Corollary}
\theoremstyle{definition}
\newtheorem{definition}[theorem]{Definition}
\newtheorem{example}[theorem]{Example}

\theoremstyle{remark}
\newtheorem{remark}[theorem]{Remark}

\DeclareMathOperator{\inp}{in}
\DeclareMathOperator{\out}{out}

\usepackage{comment}
\usepackage{subcaption}
\usepackage{bbm}

\newcommand{\pdeg}{\operatorname{pdeg}}
\newcommand{\Id}{\mathrm{Id}}
\newcommand{\NB}{\mathcal{NB}}
\newcommand{\cV}{\mathcal{V}}
\newcommand{\cE}{\mathcal{E}}
\newcommand{\cA}{\mathcal{A}}
\newcommand{\cR}{\mathcal{R}}
\newcommand{\cD}{\mathcal{D}}
\newcommand{\cL}{\mathcal{L}}
\newcommand{\cG}{\mathcal{G}}
\newcommand{\R}{\mathbb{R}}
\newcommand{\C}{\mathbb{C}}

\usepackage[normalem]{ulem}

\usepackage{tikz}
\usetikzlibrary{arrows.meta,calc,positioning,topaths}

\definecolor{color1}{RGB}{0,150,150}
\definecolor{color2}{RGB}{240,110,80}

\tikzset{
	vtx/.style={circle,fill=black,inner sep=0pt,minimum size=4.6pt},
	nbv/.style={circle,fill=color1,inner sep=0pt,minimum size=4.6pt},
	ed/.style={line width=0.7pt,black},
	oe/.style={-{Triangle[length=2.4mm,width=2mm]},line width=0.7pt,color1},
	nba/.style={-{Triangle[length=2.4mm,width=2mm]},line width=0.7pt,color2},
	come/.style={-{Triangle[length=2.4mm,width=2mm]},line width=1pt,color1},
	go/.style={-{Triangle[length=2.4mm,width=2mm]},line width=1pt,color2},
	nogo/.style={-{Triangle[length=2.4mm,width=2mm]},line width=0.7pt,black!42,dashed},
	el/.style={font=\small,color1,inner sep=1pt},
	lab/.style={font=\small,inner sep=1pt},
	qm/.style={font=\small,black!55,fill=white,inner sep=1.2pt},
	xx/.style={font=\small,black!55,fill=white,inner sep=0.5pt},
}

\definecolor{color1}{HTML}{5FA8A0}   
\definecolor{color2}{HTML}{B074A3}   
\definecolor{nodec}{HTML}{3A2A3E}    
\definecolor{edgec}{HTML}{A9A3AA}    
\definecolor{nogoc}{HTML}{A6A6A6}    

\tikzset{
	ed/.style   ={line width=0.9pt, edgec, line cap=round},
	vtx/.style  ={circle, fill=nodec, inner sep=0pt, minimum size=6.5pt},
	lab/.style  ={font=\normalsize},
	qm/.style   ={font=\Large, color2},
	el/.style   ={font=\small},
	xx/.style   ={font=\Large, nogoc},
	come/.style ={-{Triangle[length=8pt,width=7pt]}, line width=1.8pt, color1, line cap=round},
	go/.style   ={-{Triangle[length=8pt,width=7pt]}, line width=1.8pt, color2, line cap=round},
	nogo/.style ={-{Triangle[length=7pt,width=6pt]}, line width=1.4pt, nogoc, dash pattern=on 4pt off 3pt},
}

\newcommand{\cyclist}[3][1]{%
	\begin{scope}[shift={#2}, rotate=#3, scale=#1, line cap=round, line join=round]
		\draw[line width=1.5pt, color1!80!nodec] (-0.40,0.27) circle (0.25) (0.40,0.27) circle (0.25);
		\draw[line width=1.3pt, color1!80!nodec]
		(-0.40,0.27) -- (-0.02,0.27) -- (-0.14,0.60) -- cycle
		(-0.14,0.60) -- (0.28,0.60) -- (-0.02,0.27)
		(0.28,0.60) -- (0.40,0.27)
		(0.24,0.70) -- (0.30,0.60);
		\draw[line width=2.0pt, color1]
		(-0.14,0.66) -- (0.06,1.00)               
		(0.06,1.00) -- (0.27,0.76)                
		(-0.14,0.66) -- (0.10,0.52) -- (0.00,0.20); 
		\fill[color1] (0.10,1.16) circle (0.12);     
\end{scope}}

\tikzset{
	hedge/.style ={fill=#1, draw=#1, line width=24pt, line join=round},
	hedge1/.style={hedge=color1},
	hedge2/.style={hedge=color2},
}
\newcommand{\hedgethree}[4][hedge1]{%
	\begin{scope}[transparency group, opacity=0.17]
		\path[#1] (#2)--(#3)--(#4)--cycle;
	\end{scope}
	\draw[ed, line join=round] (#2)--(#3)--(#4)--cycle;}

\tikzset{
	lam3/.style={fill=color2!7, draw=color2!45, line width=0.8pt, rounded corners=9pt},
	lam2/.style={fill=color2!15, draw=color2!55, line width=0.8pt, rounded corners=7pt},
	lam1/.style={fill=color1!28, draw=color1!70, line width=0.8pt, rounded corners=5pt},
	tr/.style={black!35, line width=0.6pt},
	in1/.style={circle, fill=color1!70, draw=color1!90!black, inner sep=0pt, minimum size=6pt},
	in2/.style={circle, fill=color2!60, draw=color2!80!black, inner sep=0pt, minimum size=6pt},
	in3/.style={circle, fill=color2!30, draw=color2!60!black, inner sep=0pt, minimum size=6pt},
}

\title{Non-backtracking Laplacian for hypergraphs}

\usepackage[affil-it]{authblk}
\usepackage{microtype}

\author[1]{Eleonora Andreotti\thanks{eandreotti@fbk.eu}}
\author[2]{Raffaella Mulas\thanks{r.mulas@vu.nl}}
\affil[1]{Fondazione Bruno Kessler, Trento, Italy}
\affil[2]{Vrije Universiteit Amsterdam, Amsterdam, The Netherlands}

\date{}

\begin{document}
	\maketitle
	\begin{abstract}
		We extend the non-backtracking Laplacian for graphs, introduced by Jost, Mulas and Torres, to hypergraphs. The operator is defined on the permutations of the edges and is normalized by the number of edges through which a non-backtracking walk can continue. We show that it satisfies a symmetry that generalizes the graph case, and that it can be obtained from the random walk Laplacian of a directed graph by rescaling if and only if all edges have the same cardinality. Moreover, we show that, up to copies of the eigenvalue $1$, its spectrum can be computed from a much smaller matrix indexed by the incidences of the hypergraph, and we use this reduction to bound the remaining eigenvalues away from $1$.
	\end{abstract}
	
	\section{Introduction}

	Spectral graph theory studies the structural properties of a graph that are encoded in the eigenvalues of its associated operators, such as the adjacency matrix and the various Laplacians \cite{chung, brouwer,BookJMZ}. Among these, the non-backtracking matrix is of theoretical interest \cite{glover2021spectral,torres-1,Leo2007} and has proved remarkably effective in applications ranging from community detection \cite{bordenave2015,krzakala2013,newman2013,bolla2025clustering,bolla2025non} to node centrality \cite{arrigo2018,linzhang2019,martin2014,torres2021} and percolation \cite{shrestha2015,castellano2018, morone2015}. Moreover, in \cite{NB-Laplacian}, Jost, Mulas and Torres introduced the non-backtracking Laplacian and showed, computationally for graphs with a small number of nodes, that cospectral pairs are rarer for this operator than for the adjacency matrix, the Kirchhoff Laplacian, the normalized Laplacian and the non-backtracking matrix. Further spectral properties of this operator were established in \cite{MZZ}, and a graph coloring number arising from that work was studied in \cite{maximal-colouring}.

	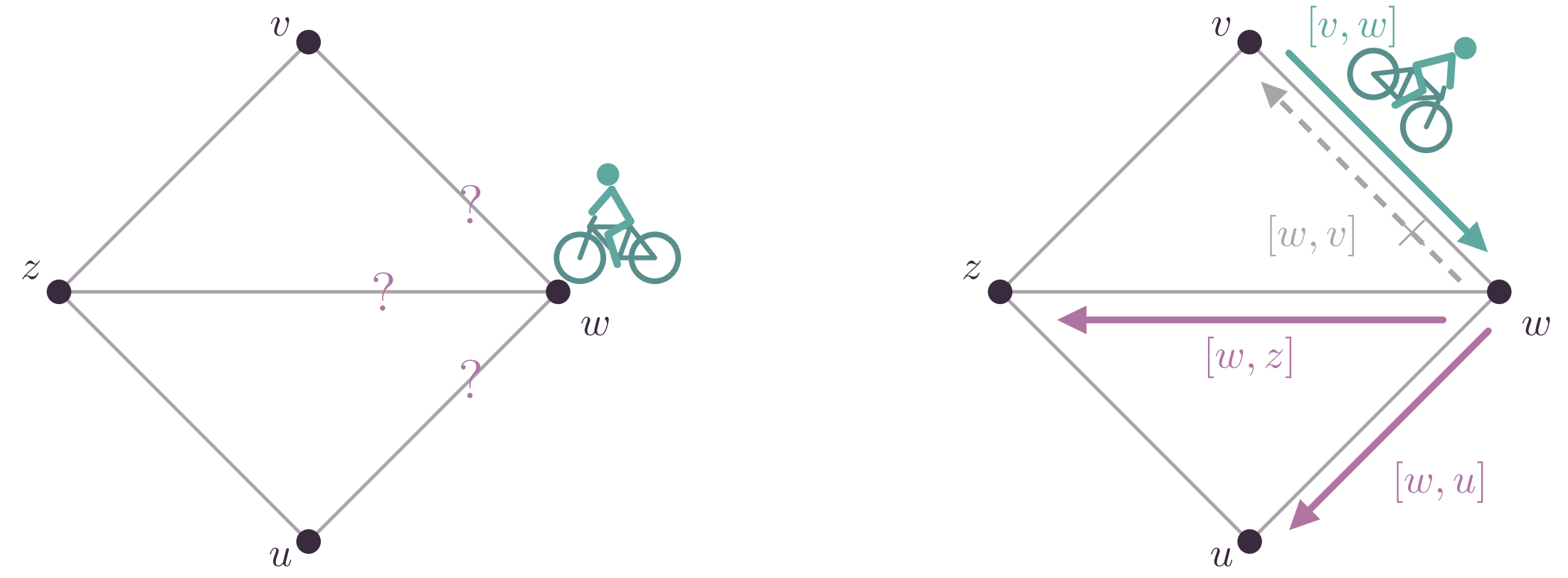
\begin{figure}[h]\centering
		\begin{subfigure}[b]{0.47\textwidth}\centering
			\begin{tikzpicture}[scale=0.58]
			\path (-5.2,-5.0) rectangle (6.2,5.2);
			\coordinate (V) at (0,4);  \coordinate (Z) at (-4,0);
			\coordinate (W) at (4,0);  \coordinate (U) at (0,-4);
			\draw[ed] (V)--(W)--(U)--(Z)--(V) (Z)--(W);
			\foreach \p in {V,Z,U}{\node[vtx] at (\p) {};}
			\node[lab,nodec] at (-0.45, 4.25) {$v$};
			\node[lab,nodec] at (-4.45, 0.35) {$z$};
			\node[lab,nodec] at ( 4.60,-0.55) {$w$};
			\node[lab,nodec] at (-0.45,-4.25) {$u$};
			\node[qm] at ( 2.60, 1.40) {$?$};
			\node[qm] at ( 1.20, 0.00) {$?$};
			\node[qm] at ( 2.60,-1.40) {$?$};
			\node[vtx] at (W) {};
			\begin{scope}[shift={((4.95,0.14)}, xscale=-1]
			\cyclist[1.5]{(0,0)}{0}
			\end{scope}
			\end{tikzpicture}
			\caption{The cyclist is riding on the vertices, and does not remember where they came from.}\label{fig:walk-a}
		\end{subfigure}\hfill
		\begin{subfigure}[b]{0.47\textwidth}\centering
			\begin{tikzpicture}[scale=0.58]
			\path (-5.2,-5.0) rectangle (6.2,5.2);
			\coordinate (V) at (0,4);  \coordinate (Z) at (-4,0);
			\coordinate (W) at (4,0);  \coordinate (U) at (0,-4);
			\draw[ed] (V)--(W)--(U)--(Z)--(V) (Z)--(W);
			\draw[come] (0.626,3.826) -- (3.826,0.626);      
			\draw[go]   (3.826,-0.626) -- (0.626,-3.826);    
			\draw[go]   (3.1,-0.45) -- (-3.1,-0.45);         
			\draw[nogo] (3.374,0.174) -- (0.174,3.374);      
			\node[xx] at (2.600,0.948) {$\times$};
			\foreach \p in {V,Z,W,U}{\node[vtx] at (\p) {};}
			\node[lab,nodec] at (-0.45, 4.25) {$v$};
			\node[lab,nodec] at (-4.45, 0.35) {$z$};
			\node[lab,nodec] at ( 4.60,-0.55) {$w$};
			\node[lab,nodec] at (-0.45,-4.25) {$u$};
			\node[lab,color1]   at ( 1.65, 4.25) {$[v,w]$};
			\node[lab,color2]   at ( 3.05,-3.05) {$[w,u]$};
			\node[lab,color2]   at ( 0.00,-1.05) {$[w,z]$};
			\node[lab,nogoc] at ( 1.00, 0.90) {$[w,v]$};
			\cyclist[1.5]{(2.12,2.78)}{-45}
			\end{tikzpicture}
			\caption{The cyclist is riding on the oriented edges, and does not need to remember.}\label{fig:walk-b}
		\end{subfigure}
		\caption{A non-backtracking walk is not Markovian as a process on the vertices, but it becomes Markovian as a process on the oriented edges. In (a) the cyclist only knows that they are at the vertex $w$, and this is not enough to
			decide along which of the three edges at $w$ they may continue. In (b) their position is the oriented edge $[v,w]$, which already records the vertex they came from: they may continue along $[w,z]$ or along $[w,u]$ (in purple), but not
			back along $[w,v]$ (dashed).}
		\label{fig:nbwalk}
	\end{figure}
	
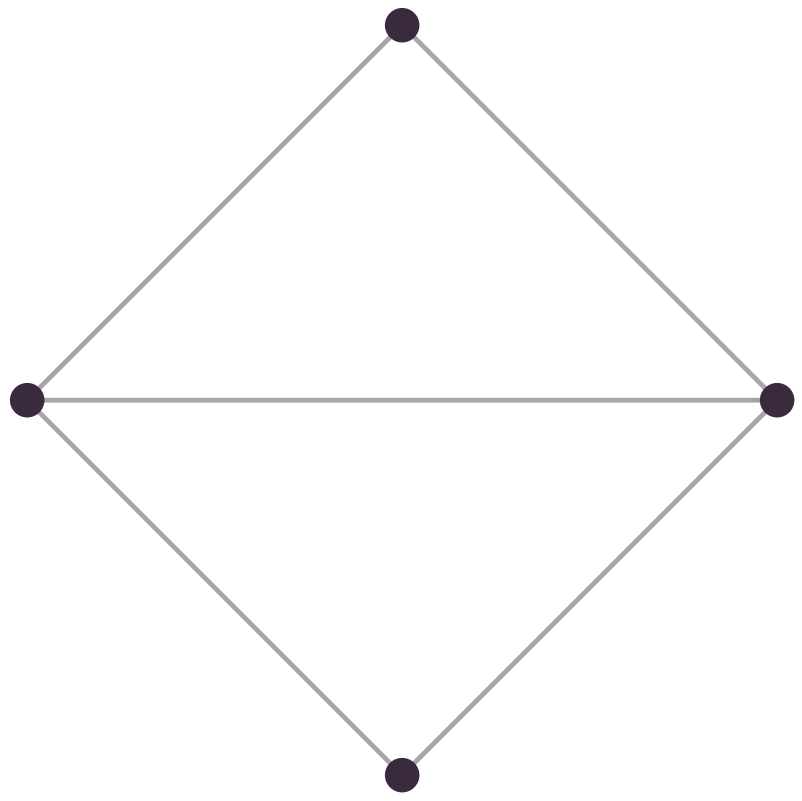
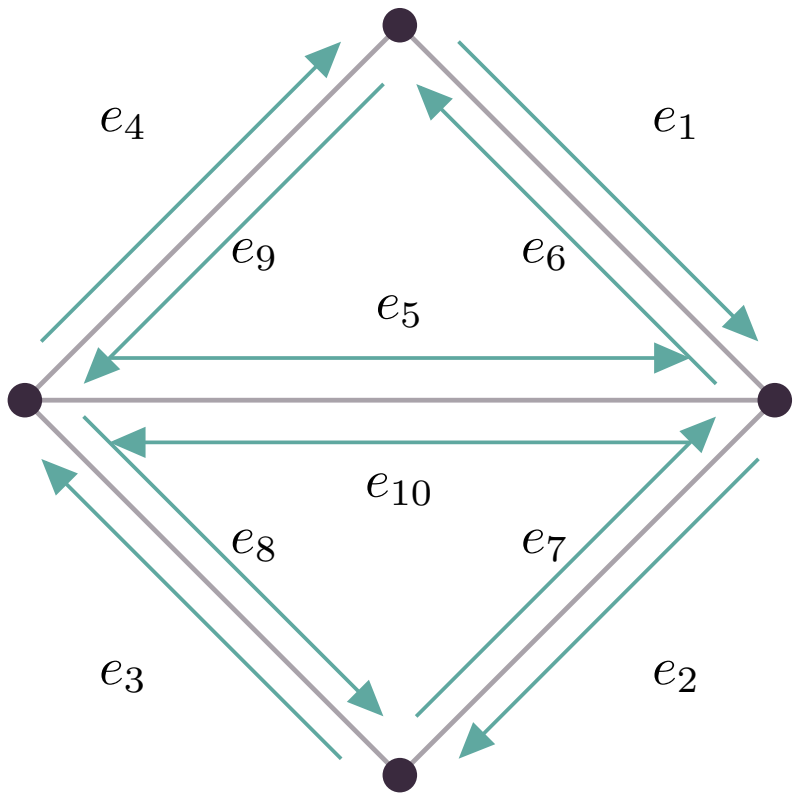
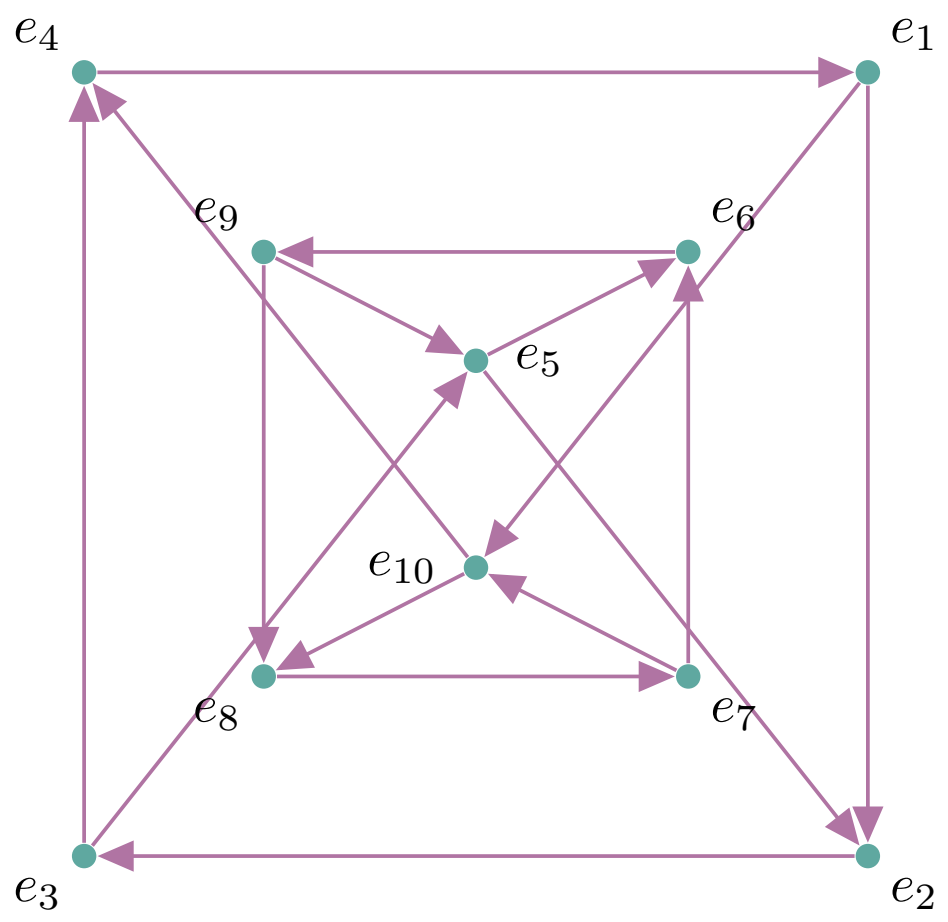
\begin{figure}[h]\centering
		\begin{subfigure}[b]{0.42\textwidth}\centering
			\begin{tikzpicture}[scale=0.62]
			\path (-4.4,-4.4) rectangle (4.4,4.4);  
			\coordinate (T) at (0,4);  \coordinate (L) at (-4,0);
			\coordinate (R) at (4,0);  \coordinate (B) at (0,-4);
			\draw[ed] (T)--(R)--(B)--(L)--(T) (L)--(R);
			\foreach \p in {T,L,R,B}{\node[vtx] at (\p) {};}
			\end{tikzpicture}
			\caption{A graph $G$ on $4$ nodes and $5$ edges.}\label{fig:nb-a}
		\end{subfigure}\hfill
		\begin{subfigure}[b]{0.42\textwidth}\centering
			\begin{tikzpicture}[scale=0.62]
			\path (-4.4,-4.4) rectangle (4.4,4.4);   
			\coordinate (T) at (0,4);  \coordinate (L) at (-4,0);
			\coordinate (R) at (4,0);  \coordinate (B) at (0,-4);
			\draw[ed] (T)--(R)--(B)--(L)--(T) (L)--(R);
			\draw[oe] (0.626,3.826) -- (3.826,0.626);    
			\draw[oe] (3.374,0.174) -- (0.174,3.374);    
			\draw[oe] (-3.826,0.626) -- (-0.626,3.826);  
			\draw[oe] (-0.174,3.374) -- (-3.374,0.174);  
			\draw[oe] (-0.626,-3.826) -- (-3.826,-0.626);
			\draw[oe] (-3.374,-0.174) -- (-0.174,-3.374);
			\draw[oe] (3.826,-0.626) -- (0.626,-3.826);  
			\draw[oe] (0.174,-3.374) -- (3.374,-0.174);  
			\draw[oe] (-3.1,0.45) -- (3.1,0.45);           
			\draw[oe] (3.1,-0.45) -- (-3.1,-0.45);         
			\foreach \p in {T,L,R,B}{\node[vtx] at (\p) {};}
			\node[el] at ( 2.95, 2.95) {$e_1$};
			\node[el] at ( 2.95,-2.95) {$e_2$};
			\node[el] at (-2.95,-2.95) {$e_3$};
			\node[el] at (-2.95, 2.95) {$e_4$};
			\node[el] at ( 1.55, 1.55) {$e_6$};
			\node[el] at ( 1.55,-1.55) {$e_7$};
			\node[el] at (-1.55,-1.55) {$e_8$};
			\node[el] at (-1.55, 1.55) {$e_9$};
			\node[el] at ( 0.00, 0.95) {$e_5$};
			\node[el] at ( 0.00,-0.95) {$e_{10}$};
			\end{tikzpicture}
			\caption{The $10$ oriented edges of $G$.}\label{fig:nb-b}
		\end{subfigure}
		
		\vspace{1.2em}
		
		\begin{subfigure}[b]{0.62\textwidth}\centering
			\begin{tikzpicture}[scale=0.72]
			\node[nbv] (x1) at ( 3.6, 3.6) {};
			\node[nbv] (x2) at ( 3.6,-3.6) {};
			\node[nbv] (x3) at (-3.6,-3.6) {};
			\node[nbv] (x4) at (-3.6, 3.6) {};
			\node[nbv] (x6) at ( 1.95, 1.95) {};
			\node[nbv] (x9) at (-1.95, 1.95) {};
			\node[nbv] (x8) at (-1.95,-1.95) {};
			\node[nbv] (x7) at ( 1.95,-1.95) {};
			\node[nbv] (x5)  at (0, 0.95) {};
			\node[nbv] (x10) at (0,-0.95) {};
			\draw[nba] (x4)--(x1);  \draw[nba] (x1)--(x2);
			\draw[nba] (x2)--(x3);  \draw[nba] (x3)--(x4);
			\draw[nba] (x6)--(x9);  \draw[nba] (x9)--(x8);
			\draw[nba] (x8)--(x7);  \draw[nba] (x7)--(x6);
			\draw[nba] (x9)--(x5);  \draw[nba] (x3)--(x5);
			\draw[nba] (x5)--(x6);  \draw[nba] (x5)--(x2);
			\draw[nba] (x1)--(x10); \draw[nba] (x7)--(x10);
			\draw[nba] (x10)--(x4); \draw[nba] (x10)--(x8);
			\node[el,above right=-2pt of x1] {$e_1$};
			\node[el,below right=-2pt of x2] {$e_2$};
			\node[el,below left=-2pt of x3]  {$e_3$};
			\node[el,above left=-2pt of x4]  {$e_4$};
			\node[el,above right=-2pt of x6] {$e_6$};
			\node[el,above left=-2pt of x9]  {$e_9$};
			\node[el,below left=-2pt of x8]  {$e_8$};
			\node[el,below right=-2pt of x7] {$e_7$};
			\node[el,right=1pt of x5]  {$e_5$};
			\node[el,left=1pt of x10]  {$e_{10}$};
			\end{tikzpicture}
			\caption{The non-backtracking graph $\mathcal{NB}(G)$.}\label{fig:nb-c}
		\end{subfigure}
		
		\caption{From a graph to its non-backtracking graph: the nodes of
			$\mathcal{NB}(G)$ are the oriented edges of $G$, with a directed edge from
			$e_i$ to $e_j$ when the output of $e_i$ is the input of $e_j$ and
			$e_j\neq e_i^{-1}$.}
		\label{fig:nbconstruction}
	\end{figure}

	The idea behind non-backtracking operators for graphs is the following. Fix a simple graph $G=(V,E)$ on $N$ nodes and $M$ edges, with minimum degree $\geq 2$. Choosing an \emph{orientation} for an edge means letting one of its endpoints be its \emph{input} and the other its \emph{output}, and we write $e=[v,w]$ for the oriented edge with input $\inp(e)=v$ and output $\out(e)=w$. Moreover, we let $e^{-1}:=[w,v]$ be the inverse of $e$, and we let $e_1,\ldots,e_{2M}$ denote the $2M$ oriented edges of $G$. A \emph{non-backtracking random walk} on $G$ is the Markov process on the oriented edges in which the walker at $e_i$ moves to each $e_j$ with $\out(e_i)=\inp(e_j)$ and $e_j\neq e_i^{-1}$ with equal probability $1/(\deg(\out(e_i))-1)$, and the minimum degree assumption guarantees that at least one such $e_j$ exists. Seen as a process on $V$ the same walk is not Markovian, since the steps available at a vertex depend on how the walker arrived there, which is why we look at the oriented edges (Figure \ref{fig:nbwalk}). The \emph{non-backtracking graph} of $G$ (Figure \ref{fig:nbconstruction}), denoted $\mathcal{NB}(G)$, was introduced in 1989 by Hashimoto \cite{hashimoto1989zeta} as the directed graph whose vertices are the $2M$ oriented edges of $G$, having a directed edge from $e_i$ to $e_j$ precisely when the above transition is admissible, that is, when $\out(e_i)=\inp(e_j)$ and $e_j\neq e_i^{-1}$. The transpose of the adjacency matrix of $\mathcal{NB}(G)$ is the \emph{non-backtracking matrix} of $G$, while the random walk Laplacian of $\mathcal{NB}(G)$ is the \emph{non-backtracking Laplacian} of $G$ \cite{NB-Laplacian}.
	
	Hypergraphs, in which an edge may join any number of vertices, are the natural model for higher-order interactions, and their spectral theory has been developed for several operators; see, for instance, \cite{BANERJEE202182,BookJMZ} and the references therein. Non-backtracking techniques for hypergraphs have appeared as well: operators defined on the \emph{incidences} of a hypergraph, that is, the edge--vertex pairs $(e,v)$ with $v\in e$, go back to the zeta function of
	Storm \cite{Storm}, and underlie spectral clustering and detectability methods \cite{ACKZ,CEH,SZ} as well as the theory of percolation on hypergraphs \cite{bianconi2024}. These operators, however, are all of adjacency type, and to the best of our knowledge no Laplacian analogue has been proposed so far.
	
	The closest object appears in the literature on hypergraph sampling: Nakajima, Kodakari and Aida \cite{nakajima2025} consider a non-backtracking walk that avoids both the previously visited node and the previously traversed edge, and derive its stochastic transition matrix on the incidences, together with its stationary distribution, with the aim of obtaining unbiased estimators of structural properties of large hypergraphs. Such a transition matrix does induce a Laplacian; being row-stochastic by construction, however, it is the random walk Laplacian of a weighted directed graph on the incidences, and its spectral theory is therefore that of a directed graph. This is an instance of a general phenomenon: transition probabilities are pairwise relations between states, and every hypergraph random walk Laplacian coincides with the random walk Laplacian of a weighted graph \cite[Theorem 4.5]{Mulas2022RandomWalks}. An operator that is not reducible to a graph Laplacian must therefore fail to be row-stochastic, and the one we introduce below does. When and why higher-order models resist reduction to graph descriptions is currently under active discussion \cite{torres2026,moutuou2026,rommens2026}.
	
	In this paper, we propose and study a generalization of the non-backtracking Laplacian for graphs in \cite{NB-Laplacian,MZZ} to the case of hypergraphs. In line with the above discussion, a hypergraph Laplacian that cannot be reduced to the graph case must give up row stochasticity, and the off-diagonal entries of our operator are therefore not always probabilities.

	\subsection*{Structure of the paper}
	
	In Section~\ref{sec:def} we introduce the non-backtracking graph, matrix and Laplacian of a hypergraph, which reduce to those of~\cite{NB-Laplacian} for graphs, and in Section~\ref{sec:NBgraph} we compute the degrees and the size of the non-backtracking graph. In Section~\ref{sec:basic} we establish the first spectral properties of the non-backtracking Laplacian, including a symmetry that generalizes that of the graph case, and we show that it is an affine function of the random walk Laplacian of the non-backtracking graph if and only if the hypergraph is uniform. Finally, in Section~\ref{sec:incidences}, we study the eigenvalue $1$, we show that the remaining eigenvalues can be computed from a smaller matrix, indexed by the incidences of the hypergraph, and we bound these eigenvalues away from $1$.

	\section{Non-backtracking operators for hypergraphs}\label{sec:def}
	
	Throughout the paper we fix a hypergraph $G=(V,E)$ on $N$ nodes and $M$ edges. We let $e_1,\ldots,e_M$ denote the edges of $G$ and, given $i\in\{1,\ldots,M\}$, we let $c_i:=|e_i|$ denote the cardinality of $e_i$. We assume that $G$ has no repeated edges, that is, \ $e_i\neq e_j$ for $i\neq j$. 
	
	The \emph{degree} of a vertex $v\in V$, denoted $\deg v$, is the number of edges in which it is contained. We assume, from here on, that $G$ has minimum vertex degree $\delta\geq 2$
	and minimum edge cardinality $\kappa\geq 2$. We also let
	$$K:=\max_i c_i \quad \text{and} \quad \Delta:=\max_{v\in V}\deg v.$$
	
	Moreover, given $i\in\{1,\ldots,M\}$, we let
	\[
	e_i^{(1)},\ldots,e_i^{(c_i!)}
	\]
	denote the $c_i!$ permutations of the vertices in $e_i$, listed with a fixed
	order.
	
	\begin{definition}\label{def:inout}
		Given a permutation $e_i^{(\alpha)}$ of $e_i$, we say that its first entry is
		its \emph{input}, and we denote it by $\inp(e_i^{(\alpha)})$; we say that its
		last entry is its \emph{output}, and we denote it by $\out(e_i^{(\alpha)})$.
		Moreover, we let
		\[
		d\bigl(e_i^{(\alpha)}\bigr) := \deg\bigl(\out(e_i^{(\alpha)})\bigr)-1 .
		\]
	\end{definition}
	
	\begin{example}\label{ex:threeedge}
		If $e_i=\{v_1,v_2,v_3\}$, then $c_i!=3!=6$, and we can let
		\begin{alignat*}{2}
			e_i^{(1)} &:= (v_1,v_2,v_3), \qquad & e_i^{(2)} &:= (v_1,v_3,v_2),\\
			e_i^{(3)} &:= (v_2,v_1,v_3), \qquad & e_i^{(4)} &:= (v_2,v_3,v_1),\\
			e_i^{(5)} &:= (v_3,v_1,v_2), \qquad & e_i^{(6)} &:= (v_3,v_2,v_1).
		\end{alignat*}
		With this fixed order of the permutations of $e_i$, we have that
		\[
		v_1=\inp(e_i^{(1)})=\inp(e_i^{(2)})=\out(e_i^{(4)})=\out(e_i^{(6)}),
		\]
		and
		\[
		d(e_i^{(4)})=\deg\bigl(\out(e_i^{(4)})\bigr)-1=\deg(v_1)-1 .
		\]
	\end{example}
	
	\begin{remark}[Permutations as orientations]\label{rmk:graph-case}
		The permutations of the edges generalize the orientations of graph edges. In fact, if $G$ is a graph, i.e.\ $c_i=2$ for all $i$, then each edge $e_i=\{v,w\}$ has exactly two permutations, $(v,w)$ and $(w,v)$, and these
		are precisely its two orientations $[v,w]$ and $[v,w]^{-1}=[w,v]$ in the sense of \cite{NB-Laplacian}. In particular, the input and the output of Definition~\ref{def:inout} are those of the corresponding oriented edge, and
		$d\bigl(e_i^{(\alpha)}\bigr)$ is the number of edges along which a non-backtracking walk that reaches $\out\bigl(e_i^{(\alpha)}\bigr)$ through
		$e_i$ can continue.  When $c_i\geq 3$, on the other hand, $e_i$ has more than two permutations, and when $c_i\geq 4$, a permutation is no longer determined by its input and output, since it also fixes the order in which the remaining $c_i-2$ vertices appear.
		
		This generalization is not the only possible one, and several
		non-equivalent notions of orientation for hypergraph edges exist in the literature: an oriented edge may be defined as a pair of disjoint vertex sets, its tail and its head~\cite{GalloLongoPallottinoNguyen}, as an ordered
		partition of its vertices~\cite{GaleanaSanchezManrique}, or as a list of its vertices in a fixed order~\cite{AustinTao}. All of them reduce to the orientations of a graph edge when $c_i=2$. We adopt the last notion: the permutations $e_i^{(1)},\dots,e_i^{(c_i!)}$ are precisely the orientations of $e_i$ in the sense of~\cite{AustinTao}, and their first and last entries provide the input and the output needed by the
		non-backtracking walk.
	\end{remark}

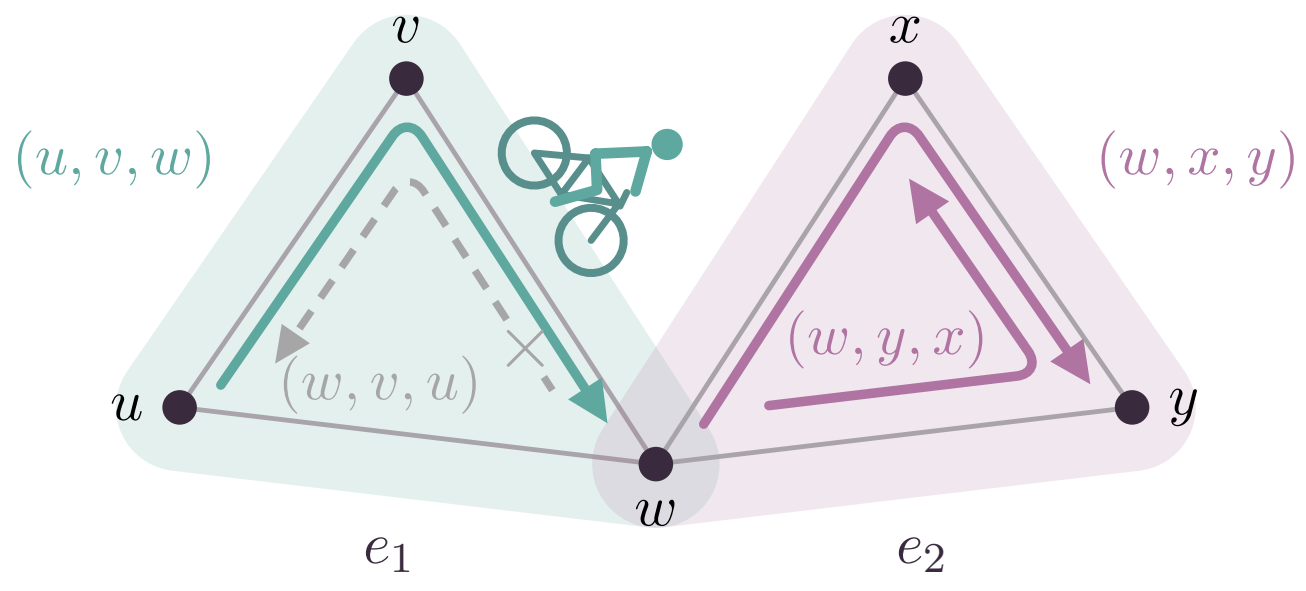
\begin{figure}[h]\centering
		\begin{tikzpicture}[scale=0.75, rc/.style={rounded corners=5pt}]
		\coordinate (U) at (-4.2,0.5);  \coordinate (V) at (-2.2,3.4);
		\coordinate (W) at (0,0);
		\coordinate (X) at (2.2,3.4);   \coordinate (Y) at (4.2,0.5);
		\coordinate (C1) at (barycentric cs:U=1,V=1,W=1);
		\coordinate (C2) at (barycentric cs:W=1,X=1,Y=1);
		\hedgethree[hedge1]{U}{W}{V}
		\hedgethree[hedge2]{W}{Y}{X}
		\draw[come,rc,shorten <=2pt,shorten >=4pt]
		($(U)!0.15!(C1)$) -- ($(V)!0.15!(C1)$) -- ($(W)!0.15!(C1)$);
		\coordinate (D0) at ($(W)!0.38!(C1)$);
		\coordinate (D1) at ($(V)!0.38!(C1)$);
		\coordinate (D2) at ($(U)!0.38!(C1)$);
		\draw[nogo,-,shorten <=4pt] (D0) -- ($(D0)!0.13!(D1)$);
		\draw[nogo,rc,shorten >=2pt] ($(D0)!0.37!(D1)$) -- (D1) -- (D2);
		\node[xx,fill=none] at ($(D0)!0.25!(D1)$) {$\times$};
		\draw[go,rc,shorten <=4pt,shorten >=2pt]
		($(W)!0.15!(C2)$) -- ($(X)!0.15!(C2)$) -- ($(Y)!0.15!(C2)$);
		\draw[go,rc,shorten <=4pt,shorten >=2pt]
		($(W)!0.38!(C2)$) -- ($(Y)!0.38!(C2)$) -- ($(X)!0.38!(C2)$);
		\foreach \p in {U,V,W,X,Y}{\node[vtx] at (\p) {};}
		\node[lab,left=3pt]  at (U) {$u$};
		\node[lab,above=3pt] at (V) {$v$};
		\node[lab,below=3pt] at (W) {$w$};
		\node[lab,above=3pt] at (X) {$x$};
		\node[lab,right=3pt] at (Y) {$y$};
		\node[lab,nodec] at ($(U)!0.5!(W)+(-0.25,-1.05)$) {$e_1$};
		\node[lab,nodec] at ($(Y)!0.5!(W)+( 0.25,-1.05)$) {$e_2$};
		\node[lab,color1,above left]  at ($(U)!0.5!(V)+(-0.5,0.3)$) {$(u,v,w)$};
		\node[lab,nogoc] at ($(C1)+(-0.3,-0.6)$) {$(w,v,u)$};
		\node[lab,color2,above right] at ($(X)!0.5!(Y)+(0.5,0.3)$)  {$(w,x,y)$};
		\node[lab,color2]             at ($(C2)+(-0.1,-0.2)$)         {$(w,y,x)$};
		\cyclist[1.15]{($($(V)!0.15!(C1)$)!0.40!($(W)!0.15!(C1)$)+(0.36,0.26)$)}{-57}
		\end{tikzpicture}
		\caption{A step of the non-backtracking walk on a hypergraph. The cyclist's
			position is the permutation $(u,v,w)$ of $e_1$, whose output is $w$. They may
			continue along any permutation of a different edge whose input is $w$,
			namely $(w,x,y)$ and $(w,y,x)$ (in purple), but not along a permutation of
			$e_1$, such as $(w,v,u)$ (dashed).}
		\label{fig:nbwalk-hypergraph}
	\end{figure}

	Generalizing the construction for graphs, we now define the non-backtracking graph of $G$ by letting the permutations of the edges play the role of the oriented edges. More precisely, a non-backtracking walk on $G$ moves from a permutation of $e_i$ to a permutation of a \emph{different} edge $e_j$, provided that the input of the latter is the output of the former (Figure \ref{fig:nbwalk-hypergraph}).

	\begin{definition}\label{def:NBgraph}
		The \emph{non-backtracking graph} of $G$ is the directed graph
		$\cG:=\NB(G)=(\cV,\cE)$ such that:
		\begin{itemize}
			\item $\cV:=\bigl\{e_i^{(\alpha)} : i\in\{1,\ldots,M\}
			\text{ and } \alpha\in\{1,\ldots,c_i!\}\bigr\}$;
			\item There exists a directed edge $e_i^{(\alpha)}\to e_j^{(\beta)}$ in
			$\cE$ if and only if
			\[
			\out\bigl(e_i^{(\alpha)}\bigr)=\inp\bigl(e_j^{(\beta)}\bigr)
			\text{ and } i\neq j .
			\]
		\end{itemize}
	\end{definition}

	Note that $\cG$ has $|\cV|=\sum_{i=1}^{M}c_i!$ nodes, which reduces to the $2M$ oriented edges when $G$ is a graph.
	
	We now generalize the non-backtracking matrix and the non-backtracking Laplacian of a graph as follows.

	\begin{definition}\label{def:NBmatrix}
		Let $\cA$ be the adjacency matrix of $\cG$, i.e.\ the $|\cV|\times|\cV|$
		matrix with $(0,1)$--entries such that
		\[
		\cA_{e_i^{(\alpha)},\,e_j^{(\beta)}} = 1
		\iff
		\out\bigl(e_i^{(\alpha)}\bigr)=\inp\bigl(e_j^{(\beta)}\bigr)
		\;\text{ and }\; i\neq j .
		\]
		The \emph{non-backtracking matrix} of $G$ is $\cA^\top$.
	\end{definition}

	\begin{definition}\label{def:NBLaplacian}
		The \emph{non-backtracking Laplacian} of $G$ is the $|\cV|\times|\cV|$
		matrix
		\[
		\cL := \Id - \cD^{-1}\cA ,
		\]
		where $\Id$ denotes the identity matrix and
		\[
		\cD := \operatorname{diag}\Bigl(
		d\bigl(e_1^{(1)}\bigr),\ldots,d\bigl(e_1^{(c_1!)}\bigr),
		\;\ldots,\;
		d\bigl(e_M^{(1)}\bigr),\ldots,d\bigl(e_M^{(c_M!)}\bigr)\Bigr).
		\]
	\end{definition}
	
	Note that $\cD$ is invertible, since the assumption that $\delta\geq 2$ implies that
	$$d(e_i^{(\alpha)})=\deg(\out(e_i^{(\alpha)}))-1\geq 1$$ for every
	$e_i^{(\alpha)}\in\cV$. Hence, $\cL$ is well-defined.
	
	\begin{remark}[Consistency with the graph case]\label{rem:graphcase}
		If $G$ is a graph, then Definitions~\ref{def:NBgraph}, \ref{def:NBmatrix}
		and~\ref{def:NBLaplacian} recover the non-backtracking graph, matrix and
		Laplacian of \cite{NB-Laplacian}. In fact, by Remark~\ref{rmk:graph-case}, in this case the nodes of $\cG$ are the oriented edges of $G$, and two oriented edges $e,f$ with $\out(e)=\inp(f)$ belong to distinct edges if and only if $\inp(e)\neq\out(f)$. Moreover, $d\bigl(e_i^{(\alpha)}\bigr)$ is the out-degree of $e_i^{(\alpha)}$ in $\cG$; therefore $\cL$ is the random walk Laplacian of $\cG$.
	\end{remark}

	\section{Degrees and size of the non-backtracking graph}\label{sec:NBgraph}
	
	In this section we study the first structural properties of the non-backtracking graph $\cG$.
	
	Recall from Definition~\ref{def:inout} that
	\[
	d\bigl(e_i^{(\alpha)}\bigr)=\deg\bigl(\out(e_i^{(\alpha)})\bigr)-1
	\]
	counts the edges through which a non-backtracking walk can continue after leaving $e_i$ through its output. If $G$ is a graph, this number coincides with the out-degree of $e_i^{(\alpha)}$ in $\cG$ (Remark~\ref{rem:graphcase}). We now compute the in- and out-degrees of the nodes of $\cG$, and we show that $d\bigl(e_i^{(\alpha)}\bigr)$ coincides with the out-degree $\deg^{\mathrm{out}}_{\cG}\bigl(e_i^{(\alpha)}\bigr)$ for every node of $\cG$ only in the graph case. In general, the out-degree counts the permutations, rather than the edges, through which the walk can continue. To count these permutations, we introduce the following notion. 
	
	\begin{definition}\label{def:pdeg}
		The \emph{permutation degree} of a vertex $v\in V$ is
		\[
		\pdeg v := \sum_{e_i \ni v}(c_i-1)!\,,
		\]
		i.e.\ the number of permutations $e_i^{(\alpha)}\in\cV$ whose input is $v$
		or, equivalently, the number of permutations whose output is $v$.
	\end{definition}
	
	Note that, if $G$ is a graph, then $\pdeg v=\deg v$ for every $v\in V$. 
	
	\begin{lemma}\label{lem:degrees}
		Let $e_i^{(\alpha)}\in\cV$, and let $v:=\out(e_i^{(\alpha)})$ and
		$w:=\inp(e_i^{(\alpha)})$. Then,
		\[
		\deg^{\mathrm{out}}_{\cG}\bigl(e_i^{(\alpha)}\bigr)
		= \pdeg v-(c_i-1)!
		\qquad\text{and}\qquad
		\deg^{\mathrm{in}}_{\cG}\bigl(e_i^{(\alpha)}\bigr) = \pdeg w-(c_i-1)! .
		\]
		In particular, every node of $\cG$ has out-degree at least $1$, and $\cD$ coincides with the out-degree matrix of $\cG$ if and only if $G$ is a graph.
	\end{lemma}
	\begin{proof}
		By Definition~\ref{def:NBgraph}, the out-neighbors of $e_i^{(\alpha)}$ are the permutations $e_j^{(\beta)}$ with $j\neq i$ and $\inp(e_j^{(\beta)})=v$. For each $j\neq i$ with $e_j\ni v$, there are exactly $(c_j-1)!$ permutations of $e_j$ whose input is $v$. Therefore,
		\[
		\deg^{\mathrm{out}}_{\cG}\bigl(e_i^{(\alpha)}\bigr)
		= \sum_{j\neq i,\; e_j\ni v}(c_j-1)!
		= \pdeg v-(c_i-1)! ,
		\]
		and the formula for the in-degree follows in the same way, by counting, for each $j\neq i$ with $e_j\ni w$, the $(c_j-1)!$ permutations of $e_j$ whose output is $w$. Since $\deg v\geq\delta\geq 2$, there exists at least one $j\neq i$ with $e_j\ni v$, hence the out-degree is at least $1$. Moreover, the sum above has exactly $\deg v-1=d\bigl(e_i^{(\alpha)}\bigr)$ terms, each of which is at least $1$. Hence, $\deg^{\mathrm{out}}_{\cG}\bigl(e_i^{(\alpha)}\bigr)=d\bigl(e_i^{(\alpha)}\bigr)$ if and only if $(c_j-1)!=1$, i.e.\ $c_j=2$, for every $j\neq i$ with $e_j\ni v$. Since every vertex of $e_i$ is the output of some permutation of $e_i$, letting $e_i^{(\alpha)}$ range over all of $\cV$ shows that $\cD$ coincides with the out-degree matrix of $\cG$ if and only if $c_j=2$ for every edge $e_j$ that meets some other edge. Since $\delta\geq 2$, every edge meets some other edge, hence this condition is equivalent to $c_j=2$ for all $j$, that is, to $G$ being a graph.
	\end{proof}

	\begin{remark}\label{rem:notprob}
		Fix $e_i^{(\alpha)}\in\cV$ and let $v:=\out(e_i^{(\alpha)})$. The row of $\cA$ indexed by $e_i^{(\alpha)}$ has an entry equal to $1$ for each out-neighbor of $e_i^{(\alpha)}$ in $\cG$, and all other entries equal to $0$, therefore it sums to $\deg^{\mathrm{out}}_{\cG}\bigl(e_i^{(\alpha)}\bigr)$. Multiplying by $\cD^{-1}$ divides this row by $d\bigl(e_i^{(\alpha)}\bigr)=\deg v-1$. Hence, by Lemma~\ref{lem:degrees}, the row of $\cD^{-1}\cA$ indexed by $e_i^{(\alpha)}$ sums to
		\[
		\frac{\pdeg v-(c_i-1)!}{\deg v-1}
		\;=\;
		\frac{1}{\deg v-1}\sum_{j\neq i,\;e_j\ni v}(c_j-1)!\,.
		\]
		The sum on the right has exactly $\deg v-1$ terms, one for each edge through which the walk can continue after leaving $e_i$ at $v$, and each term is at least $1$. Therefore, the row sum is at least $1$, with equality if and only if $c_j=2$ for every $j\neq i$ with $e_j\ni v$. As in the proof of Lemma~\ref{lem:degrees}, this condition holds for every node of $\cG$ if and only if $G$ is a graph. Hence, in contrast with the graph case, $\cD^{-1}\cA$ is row-stochastic, and $\cL$ is the random walk Laplacian of $\cG$, if and only if $G$ is a graph.
	\end{remark}
	
	We can now prove the following generalization of \cite[Lemma 3.1]{NB-Laplacian}.
	
	\begin{lemma}\label{lem:count}
		The non-backtracking graph $\cG$ has
		\[
		|\cV| = \sum_{i=1}^{M} c_i!
		\qquad\text{nodes and}\qquad
		|\cE| = \sum_{v\in V}\;
		\sum_{\substack{(i,j):\ i\neq j,\\ v\in e_i\cap e_j}}(c_i-1)!\,(c_j-1)!
		\qquad\text{edges.}
		\]
	\end{lemma}

	\begin{proof}
		The number of nodes was already observed after Definition~\ref{def:NBgraph}. Now, observe that $|\cE|$ equals the sum of the out-degrees of the nodes of $\cG$. We group the nodes of $\cG$ according to their edge and their output vertex: for each $v\in V$ and each $e_i\ni v$, there are exactly $(c_i-1)!$ permutations of $e_i$ whose output is $v$, and by Lemma~\ref{lem:degrees} each of them has out-degree $\pdeg v-(c_i-1)!$. Hence,
		\begin{align*}
			|\cE|
			&= \sum_{v\in V}\;\sum_{e_i\,\ni\, v}(c_i-1)!\,
			\bigl(\pdeg v-(c_i-1)!\bigr)
			\\ &= \sum_{v\in V}\;\sum_{e_i\,\ni\, v}(c_i-1)!
			\sum_{\substack{e_j\,\ni\, v\\ j\neq i}}(c_j-1)!
			\\ &= \sum_{v\in V}\;
			\sum_{\substack{(i,j):\ i\neq j,\\ v\in e_i\cap e_j}}(c_i-1)!\,(c_j-1)! .
		\end{align*}
	\end{proof}
	
	\begin{remark}
		If $G$ is a graph, then $(c_i-1)!=1$ for all $i$, therefore Lemma
		\ref{lem:count} gives
		\[
		|\cE|=\sum_{v\in V}\deg v\,(\deg v-1)
		=\sum_{v\in V}(\deg v)^2-2M,
		\]
		which is the same as in \cite[Lemma 3.1]{NB-Laplacian}.
	\end{remark}
	
	We illustrate the results above on two families of hypergraphs. The first is a family of uniform hypergraphs, that is, hypergraphs whose edges all have the same cardinality; if this cardinality is $\kappa$, we say that the hypergraph is \emph{$\kappa$-uniform}. Moreover, we say that a hypergraph is \emph{$r$-regular} if all its vertices have degree $r$.
	
\begin{figure}[htb]\centering
		\tikzset{
			hullrow/.style={fill=color1, fill opacity=0.20, rounded corners=0.27cm},
			hullcol/.style={fill=color2, fill opacity=0.20, rounded corners=0.27cm},
		}
		\newcommand{\latticebase}{%
			\foreach \a in {1,2,3}{
				\fill[hullrow] (1.4-0.28,-1.4*\a+0.28) rectangle (4.2+0.28,-1.4*\a-0.28);
				\node[lab,color1!80!nodec] at (0.6,-1.4*\a) {$e_{\a}$};}
			\foreach \b in {1,2,3}{
				\fill[hullcol] (1.4*\b-0.28,-1.4+0.28) rectangle (1.4*\b+0.28,-4.2-0.28);
				\node[lab,color2!80!nodec] at (1.4*\b,-5.0) {$e_{\the\numexpr3+\b}$};}
			\foreach \a in {1,2,3}{
				\draw[ed] (1.4,-1.4*\a) -- (4.2,-1.4*\a);
				\draw[ed] (1.4*\a,-1.4) -- (1.4*\a,-4.2);}}
		\begin{subfigure}[b]{0.42\textwidth}\centering
			\begin{tikzpicture}
			\path (0.2,-5.3) rectangle (5.4,-0.4);
			\latticebase
			\foreach \a in {1,2,3} \foreach \b in {1,2,3}{
				\node[vtx] at (1.4*\b,-1.4*\a) {};
				\node[lab,nodec,font=\scriptsize,anchor=north west,inner sep=0.5pt] at (1.4*\b+0.33,-1.4*\a-0.33) {$v_{\a\b}$};}
			\end{tikzpicture}
			\caption{The $3$-lattice.}\label{fig:lattice-a}
		\end{subfigure}\hfill
		\begin{subfigure}[b]{0.52\textwidth}\centering
			\begin{tikzpicture}[scale=0.7]
			\path (-1.6,-3.3) rectangle (7.9,3.6);
			\cyclist[1.25]{(0,0.40)}{0}
			\node[lab,color1] (p) at (0,0) {$(v_{11},v_{12},v_{13})$};
			\node[lab,color2,anchor=west] (q) at (4.2,2.3) {$(v_{13},v_{23},v_{33})$};
			\node[lab,color2,anchor=west] (r) at (4.2,0)   {$(v_{13},v_{33},v_{23})$};
			\node[lab,nogoc,anchor=west] (s) at (4.2,-2.3) {$(v_{13},v_{12},v_{11})$};
			\node[lab,font=\footnotesize,color2,anchor=north west] at (q.south west) {in $e_6$};
			\node[lab,font=\footnotesize,color2,anchor=north west] at (r.south west) {in $e_6$};
			\node[lab,font=\footnotesize,nogoc,anchor=north west] at (s.south west) {in $e_1$};
			\draw[go,shorten >=2pt]   (p.east) -- (q.west);
			\draw[go,shorten >=2pt]   (p.east) -- (r.west);
			\draw[nogo,shorten >=2pt] (p.east) -- (s.west);
			\node[xx] at ($(p.east)!0.55!(s.west)$) {$\times$};
			\end{tikzpicture}
			\caption{A step of the walk from a permutation of $e_1$.}\label{fig:lattice-b}
		\end{subfigure}
		\caption{In (a), the $3$-lattice, whose edges are the rows $e_1,e_2,e_3$ and the
			columns $e_4,e_5,e_6$. In (b), the cyclist's position is the permutation
			$(v_{11},v_{12},v_{13})$ of $e_1$, whose output is $v_{13}$. The only other
			edge through $v_{13}$ is $e_6$, so the cyclist may continue along its
			$(c_6-1)!=2$ permutations with input $v_{13}$ (in purple), but not along a
			permutation of $e_1$ (dashed).}
		\label{fig:lattice}
	\end{figure}
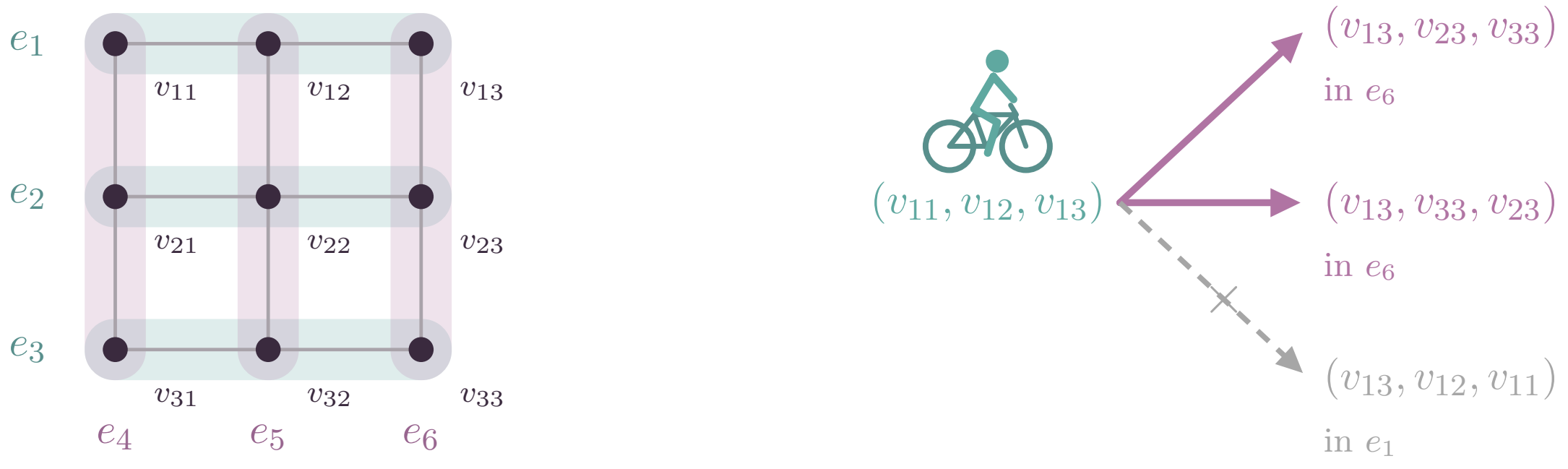

	\begin{example}\label{ex:lattice}
		Given $\ell\geq 2$, the \emph{$\ell$-lattice} is the hypergraph whose
		vertices are the $\ell^2$ cells of an $\ell\times\ell$ grid and whose $2\ell$ edges are its rows and columns~\cite{AM,A25}. It is $\ell$-uniform and $2$-regular, and for $\ell=2$, it is the cycle graph on $4$ vertices.
		
		Let $G$ be the $3$-lattice (Figure~\ref{fig:lattice-a}): its vertices are $v_{ab}$ for $a,b\in\{1,2,3\}$, while its edges are the rows $e_a=\{v_{a1},v_{a2},v_{a3}\}$ and the columns $e_{3+b}=\{v_{1b},v_{2b},v_{3b}\}$. In this case, $N=9$, $M=6$, and $c_i=3$ for all $i$. Since $\deg v=2$ for every vertex, we have $d\bigl(e_i^{(\alpha)}\bigr)=1$ for every permutation, so that $\cD=\Id$ and $\cL=\Id-\cA$. Moreover, $\pdeg v=2\cdot 2!=4$ for every $v\in V$, so that, by Lemmas~\ref{lem:degrees} and~\ref{lem:count}, every node of $\cG$ has out-degree $\pdeg v-(c_i-1)!=2$, and $\cG$ has
		\[
		|\cV|=6\cdot 3!=36 \quad\text{nodes and}\quad
		|\cE|=\sum_{v\in V} 2\cdot 2!\cdot 2!=72 \quad\text{edges,}
		\]
		since each vertex lies in exactly two edges, which give two ordered pairs $(i,j)$. The local structure of $\cG$ is shown in Figure~\ref{fig:lattice-b}: from a permutation of a row with output $v$, the walk can only continue into the unique column through $v$, along either of its $(3-1)!=2$ permutations with input $v$, and vice versa. In particular, every walk on $\cG$ alternates between permutations of rows and permutations of columns, therefore $\cG$ is bipartite.
	\end{example}

	In the lattices, all edges have the same cardinality. In the next family, the cardinalities of the edges vary.

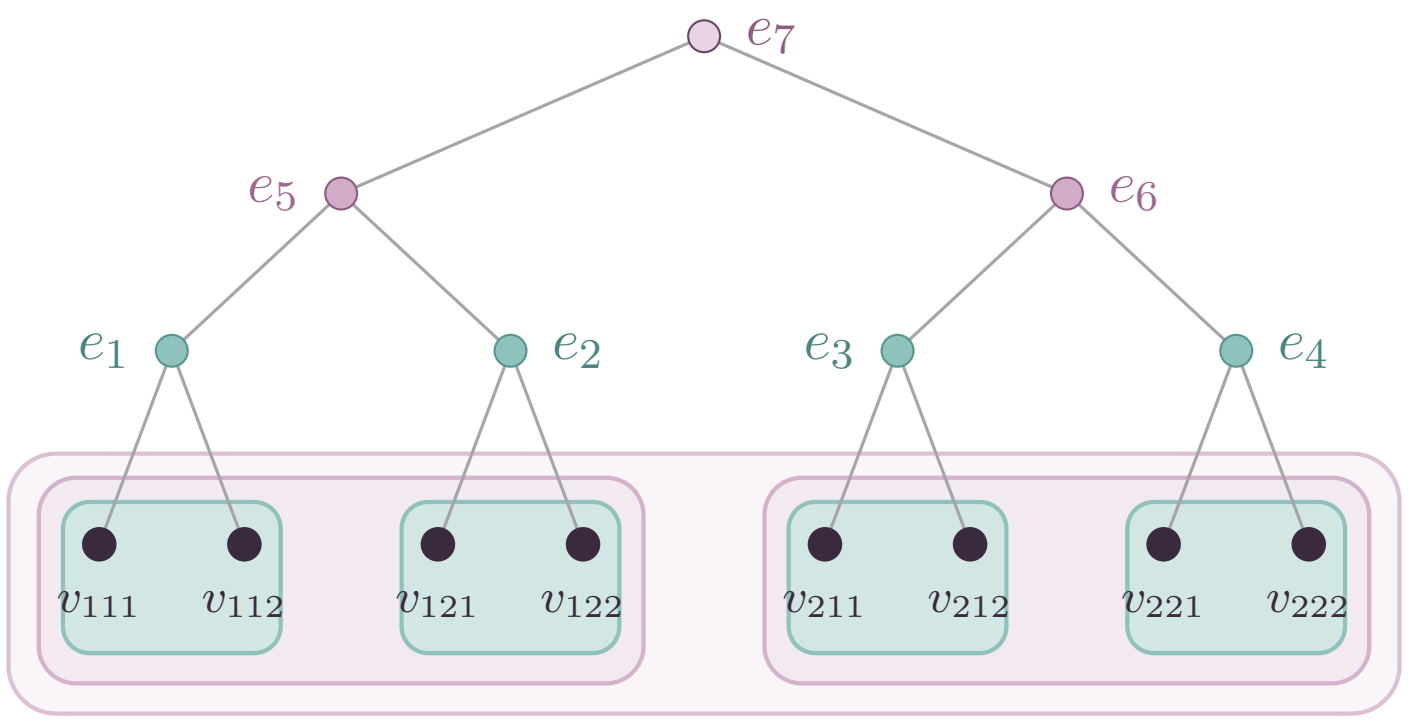
\begin{figure}[h]\centering
		\begin{tikzpicture}[scale=0.8]
		\foreach \w [count=\i from 0] in {111,112,121,122,211,212,221,222}{
			\coordinate (L\w) at ({1.2*\i+0.4*floor(\i/2)+0.4*floor(\i/4)},0);}
		\coordinate (N11) at ($(L111)!0.5!(L112)+(0,1.6)$);
		\coordinate (N12) at ($(L121)!0.5!(L122)+(0,1.6)$);
		\coordinate (N21) at ($(L211)!0.5!(L212)+(0,1.6)$);
		\coordinate (N22) at ($(L221)!0.5!(L222)+(0,1.6)$);
		\coordinate (N1)  at ($(N11)!0.5!(N12)+(0,1.3)$);
		\coordinate (N2)  at ($(N21)!0.5!(N22)+(0,1.3)$);
		\coordinate (R)   at ($(N1)!0.5!(N2)+(0,1.3)$);
		\draw[lam3] ($(L111)+(-0.75,-1.4)$) rectangle ($(L222)+(0.75,0.75)$);
		\draw[lam2] ($(L111)+(-0.5,-1.15)$) rectangle ($(L122)+(0.5,0.55)$);
		\draw[lam2] ($(L211)+(-0.5,-1.15)$) rectangle ($(L222)+(0.5,0.55)$);
		\foreach \a/\b in {111/112,121/122,211/212,221/222}{
			\draw[lam1] ($(L\a)+(-0.3,-0.9)$) rectangle ($(L\b)+(0.3,0.35)$);}
		\draw[tr] (R)--(N1) (R)--(N2)
		(N1)--(N11) (N1)--(N12) (N2)--(N21) (N2)--(N22)
		(N11)--(L111) (N11)--(L112) (N12)--(L121) (N12)--(L122)
		(N21)--(L211) (N21)--(L212) (N22)--(L221) (N22)--(L222);
		\node[in3] at (R) {};   \node[lab,color2!80!black,right=4pt] at (R) {$e_7$};
		\node[in2] at (N1) {};  \node[lab,color2!90!black,left=4pt]  at (N1) {$e_5$};
		\node[in2] at (N2) {};  \node[lab,color2!90!black,right=4pt] at (N2) {$e_6$};
		\foreach \n/\e/\s in {N11/1/left,N12/2/right,N21/3/left,N22/4/right}{
			\node[in1] at (\n) {}; \node[lab,color1!80!black,\s=4pt] at (\n) {$e_{\e}$};}
		\foreach \w in {111,112,121,122,211,212,221,222}{
			\node[vtx] at (L\w) {};
			\node[lab,nodec,font=\footnotesize,below=5pt] at (L\w) {$v_{\w}$};}
		\end{tikzpicture}
		\caption{The laminar hypergraph $\Lambda_{2,3}$ obtained from the binary tree of depth $3$. Each internal node of the tree gives the edge formed by the leaves below it, drawn in the same color: the four nodes at height $1$ give the level-$1$ edges $e_1,\dots,e_4$, the two nodes at height $2$ give the level-$2$ edges $e_5,e_6$, and the root gives the level-$3$ edge $e_7$.}
		\label{fig:laminar-tree}
	\end{figure}

	\begin{example}\label{ex:laminar} Given $k\geq 2$ and $p\geq 2$, consider the rooted tree in which every node other than the leaves has exactly $k$ children and all leaves are at distance $p$ from the root (the complete $k$-ary tree of depth $p$). The \emph{$(k,p)$-laminar hypergraph} $\Lambda_{k,p}$ has as vertices the $k^p$ leaves of the tree and, as edges, for each internal node of the tree, that is, for each node that is not a leaf, the set of leaves below it (Figure~\ref{fig:laminar-tree}). 
		
		If an internal node has height $\ell\in\{1,\dots,p\}$, that is, if its leaves are at distance $\ell$ below it, then it gives an edge of cardinality $k^\ell$, which we call a \emph{level-$\ell$ edge}. In particular, any two edges are either disjoint or one of them contains the other.
		
		Labeling the leaves by the words in $\{1,\dots,k\}^p$, the level-$\ell$ edges are the sets of words sharing the first $p-\ell$ letters. Hence, $\Lambda_{k,p}$ has $N=k^p$ nodes and $M=(k^p-1)/(k-1)$ edges, of which $k^{p-\ell}$ have cardinality $k^\ell$; since every leaf has exactly one ancestor at each height, $\deg v=p$ for all $v$, and $\Lambda_{k,p}$ is non-uniform. Since $\delta=p$ and $\kappa=k$, the hypergraph $\Lambda_{k,p}$ satisfies our assumptions.
		
		We consider $\Lambda_{2,3}$ (Figure~\ref{fig:laminar-tree}), which has $N=8$, $M=7$, and edges of cardinality $2$ (four of them, level $1$), $4$ (two of them, level $2$) and $8$ (one, level $3$). Since $\deg v=3$ for every vertex, we have that $d(e_i^{(\alpha)})=2$ for every permutation, therefore
		$\cD=2\,\Id$ and
		\[
		\cL=\Id-\tfrac12\,\cA .
		\]
		Moreover, $\pdeg v=\sum_{\ell=1}^{3}(2^{\ell}-1)!=1!+3!+7!=5047$ for every $v\in V$, hence, by Lemmas \ref{lem:degrees} and \ref{lem:count}, a permutation of a level-$\ell$ edge has out-degree
		$5047-(2^{\ell}-1)!$ in $\cG$ -- namely $5046$, $5041$ and $7$ for
		$\ell=1,2,3$ -- and $\cG$ has
		\[
		|\cV|=\sum_{\ell=1}^{3}2^{\,3-\ell}(2^{\ell})!
		=4\cdot 2!+2\cdot 4!+8!=40376
		\quad\text{nodes and}\quad\]
		\[  |\cE|=8\bigl(5047^2-(1^2+6^2+5040^2)\bigr)=564576
		\quad\text{edges.}
		\]
		
		Comparing with Example~\ref{ex:lattice} shows the effect of non-uniformity. As for the $3$-lattice, all vertices of $\Lambda_{2,3}$ have the same degree, therefore $\cD$ is a scalar matrix; yet, since $\Lambda_{2,3}$ is not uniform, the out-degrees in $\cG$ depend on the level, ranging from $7$ to $5046$. Hence, the row sums of $\cD^{-1}\cA$ equal $2523$, $2520.5$ and $3.5$ on the permutations of the level-$1$, level-$2$ and level-$3$ edges, respectively, and $\cD^{-1}\cA$ is far from being row-stochastic (Remark~\ref{rem:notprob}). 
	\end{example}

	\section{First spectral properties}\label{sec:basic}
	We now turn to the spectral properties of $\cL$.  Given a square matrix $\mathcal{O}$, we let $\sigma(\mathcal{O})$ denote its spectrum, and we let $D(z,r)\subseteq\C$ denote the closed disk of center $z$ and radius $r$.
	
	If $G$ is a graph, then by \cite[Remark 3.3]{NB-Laplacian}, $\mathcal{L}$ has $2M$ eigenvalues (counted with algebraic multiplicity) that sum to $2M$. Moreover, in this case $\sigma(\mathcal{L})\subseteq D(1,1)$ \cite{Bauer} and therefore, the real eigenvalues of $\mathcal{L}$ are contained in the interval $[0,2]$. For hypergraphs, these facts generalize as follows.

	\begin{proposition}\label{prop:gershgorin}
		The non-backtracking Laplacian $\cL$ has 
		$\sum_{i=1}^M c_i!$ eigenvalues (counted with algebraic multiplicity) that sum to 
		$\sum_{i=1}^M c_i!$. Moreover,
		\[
		\sigma(\cL)\subseteq D\bigl(1,\,(K-1)!\bigr),
		\]
		where $K$ is the largest edge cardinality. In particular, the real eigenvalues of $\cL$ are contained in the interval $\bigl[1-(K-1)!,\;1+(K-1)!\bigr]$.
	\end{proposition}
	
	\begin{proof}
		Since $\cL$ is a $|\cV|\times|\cV|$ matrix, it has $|\cV|=\sum_{i=1}^M c_i!$ eigenvalues, counted with algebraic multiplicity, and their sum is $\operatorname{tr}\cL$. Since all diagonal entries of $\cL$ equal $1$, we have that $\operatorname{tr}\cL=|\cV|$. This proves the first claim.
		
		Now, all the Gershgorin disks of $\cL$ are centered at $1$. By Lemma~\ref{lem:degrees}, the radius of the disk corresponding to the row of $\mathcal{L}$ indexed by $e_i^{(\alpha)}$, with $v:=\out(e_i^{(\alpha)})$, is
		\[
		\frac{1}{\deg v-1}\sum_{j\neq i,\;e_j \ni v}(c_j-1)!
		\;\leq\;
		\frac{(\deg v-1)\cdot(K-1)!}{\deg v-1}
		\;=\;(K-1)! ,
		\]
		since the sum has exactly $\deg v-1$ terms, each bounded above by $(K-1)!$. Hence, by the Gershgorin Disk Theorem, $$\sigma(\cL)\subseteq D\bigl(1,(K-1)!\bigr).$$
	\end{proof}

	We now generalize~\cite[Theorem 3.8]{NB-Laplacian} to hypergraphs. To this end, we give the following definitions.
	
	\begin{definition}\label{def:reversal}
		The \emph{reversal} of a permutation $\xi=(v_1,v_2,\dots,v_c)\in\cV$ is the permutation $\bar\xi:=(v_c,\dots,v_2,v_1)\in\cV$ of the same edge. We let $P$ be the $|\cV|\times|\cV|$ permutation matrix of the map $\xi\mapsto\bar\xi$, and, given $\mathbf{x},\mathbf{y}\in\C^{|\cV|}$, we define their \emph{$P$-product} as
		\[
		(\mathbf{x},\mathbf{y})_P := \langle \mathbf{x},P\mathbf{y}\rangle = \overline{\mathbf{x}}^\top P\,\mathbf{y} ,
		\]
		where $\langle \mathbf{x},\mathbf{y}\rangle:=\overline{\mathbf{x}}^\top \mathbf{y}$ is the usual complex inner product.
	\end{definition}

	With this notation, \cite[Theorem 3.8]{NB-Laplacian} generalizes as follows.

	\begin{theorem}\label{thm:PT}
		Let $\cL$ be the non-backtracking Laplacian of a hypergraph. Then,
		\begin{enumerate}
			\item $\cL^\top = P\,\cL\,P$, i.e., $\cL$ is PT-symmetric;
			\item $\cL$ is self-adjoint with respect to the $P$-product, i.e., $\cL$ is
			$P$-self-adjoint;
			\item If $(\lambda,\mathbf{x})$ is an eigenpair for $\cL$ and $\lambda\in\C$ is not real, then
			\[
			\sum_{i=1}^{M}\;\sum_{\substack{\xi\text{ permutation} \\ \text{of }e_i}}\overline{x_{\xi}}\;x_{\bar\xi}=0 .
			\]
		\end{enumerate}
	\end{theorem}
	
	\begin{proof}
		Note that, by definition, $\inp(\bar\xi)=\out(\xi)$ and $\out(\bar\xi)=\inp(\xi)$ for every $\xi\in\cV$, and that $P^\top=P$ and $P^2=\Id$, since reversing a permutation twice gives back the original one.
		\begin{enumerate}
			\item Let $\xi,\eta\in\cV$ be permutations of $e_i$ and $e_j$, respectively, so that $\bar\xi$ and $\bar\eta$ are also permutations of $e_i$ and $e_j$. We write $\cL_{\xi\eta}$ for the entry of $\cL$ in the row indexed by $\xi$ and the column indexed by $\eta$. For $\xi\neq\eta$ we have that
			\[
			\cL_{\bar\xi\bar\eta}
			= -\frac{1}{d(\bar\xi)}\,
			\mathbbm{1}\bigl[\out(\bar\xi)=\inp(\bar\eta)\bigr]\cdot\mathbbm{1}[i\neq j]
			= -\frac{1}{d(\bar\xi)}\,
			\mathbbm{1}\bigl[\inp(\xi)=\out(\eta)\bigr]\cdot\mathbbm{1}[i\neq j],
			\]
			where $\mathbbm{1}$ denotes the indicator function. When this entry is non-zero, let $u:=\inp(\xi)=\out(\eta)$. Then,
			\[
			d(\bar\xi)=\deg\bigl(\out(\bar\xi)\bigr)-1=\deg u-1=\deg\bigl(\out(\eta)\bigr)-1=d(\eta),
			\]
			so that $\cL_{\bar\xi\bar\eta}=\cL_{\eta\xi}$. On the diagonal, both $\cL_{\bar\xi\bar\xi}$ and $\cL_{\xi\xi}$ equal $1$, and $\bar\xi=\bar\eta$ if and only if $\xi=\eta$. Therefore, $$(P\cL P)_{\xi\eta}=\cL_{\bar\xi\bar\eta}=\cL_{\eta\xi}=(\cL^\top)_{\xi\eta}$$ for all $\xi,\eta\in\cV$, which proves the first claim.
			\item Since $P^\top=P$ and $P^2=\Id$, we have that
			\[
			(\mathbf{x},\cL \mathbf{y})_P=\langle \mathbf{x}, P\cL \mathbf{y}\rangle
			=\langle \mathbf{x}, \cL^\top P \mathbf{y}\rangle
			=\langle \cL \mathbf{x}, P \mathbf{y}\rangle=(\cL \mathbf{x},\mathbf{y})_P ,
			\]
			where we have used that $P\cL=\cL^\top P$, which follows from the first claim, and that $\cL$ is a real matrix. Therefore, $\cL$ is self-adjoint with respect to the $P$-product.
			\item The second claim implies that, if $(\lambda,\mathbf{x})$ is an eigenpair for $\cL$, then
			\[
			\lambda(\mathbf{x},\mathbf{x})_P=(\mathbf{x},\lambda \mathbf{x})_P=(\mathbf{x},\cL \mathbf{x})_P
			=(\cL \mathbf{x},\mathbf{x})_P=(\lambda \mathbf{x},\mathbf{x})_P=\overline{\lambda}(\mathbf{x},\mathbf{x})_P .
			\]
			Hence, if $\lambda\neq\overline{\lambda}$, i.e.\ $\lambda$ is not real, then $(\mathbf{x},\mathbf{x})_P=\overline{\mathbf{x}}^\top P\mathbf{x}=0$. Since $(P\mathbf{x})_\xi=x_{\bar\xi}$ for every $\xi\in\cV$, this is the identity in the statement.
		\end{enumerate}
	\end{proof}

	\begin{remark}\label{rem:normalization-PT}
		The proof of Theorem~\ref{thm:PT} only uses the fact that $d\bigl(e_i^{(\alpha)}\bigr)$
		depends on $\out(e_i^{(\alpha)})$ alone. This singles out $\cD$ among the
		natural ways of counting how the walk can continue after exiting $e_i$ at
		$v:=\out(e_i^{(\alpha)})$, namely:
		\begin{itemize}
			\item the number $\deg v-1$ of edges $e_j\ni v$ with $j\neq i$, which is
			our choice;
			\item the number $\sum_{e_j\ni v,\,j\neq i}(c_j-1)$ of pairs $(e_j,w)$
			with $e_j\ni v$, $j\neq i$ and $w\in e_j\setminus\{v\}$;
			\item the number $\pdeg v-(c_i-1)!$ of permutations, i.e.\ the out-degree
			of $e_i^{(\alpha)}$ in $\cG$ (Lemma~\ref{lem:degrees}).
		\end{itemize}
		The last two also depend on $i$, and therefore yield a PT-symmetric
		operator if and only if any two edges sharing a vertex have the same
		cardinality, i.e.\ if and only if every connected component of $G$ is
		uniform; in this case, the three normalizations are proportional on each
		component.
		
		On the other hand, only the last one always makes $\cD^{-1}\cA$
		row-stochastic: the first does so if and only if $G$ is a graph
		(Remark~\ref{rem:notprob}), and the second if and only if all edges have
		at most $3$ vertices.
		
		Hence, unless every connected component of $G$ is uniform, none of these
		normalizations is both PT-symmetric and stochastic: our choice preserves
		the symmetry, at the price of stochasticity.
	\end{remark}
	
	We conclude this section by showing that the uniform hypergraphs are precisely those for which $\cL$ is an affine image of the random walk Laplacian of $\cG$, so that its spectrum is obtained from that of a random walk Laplacian by an affine transformation. To this end, we let $\cD_{\cG}$ denote the out-degree matrix of $\cG$, and we let
	\[
	\cL_{\mathrm{rw}}:=\Id-\cD_{\cG}^{-1}\cA
	\]
	denote the random walk Laplacian of $\cG$. Recall that $G$ is $\kappa$-uniform if all its edges have cardinality $\kappa$.

	\begin{proposition}\label{prop:uniform}
		There exist constants $a,b\in\R$ such that
		\[
		\cL = a\,\Id + b\,\cL_{\mathrm{rw}}
		\]
		if and only if $G$ is uniform. In this case, if $G$ is $\kappa$-uniform, then $a=1-(\kappa-1)!$ and $b=(\kappa-1)!$, and
		\begin{enumerate}
			\item $\sigma(\cL)$ is the image of $\sigma(\cL_{\mathrm{rw}})$ under the map $z\mapsto 1+(\kappa-1)!\,(z-1)$;
			\item the constant functions are eigenfunctions of $\cL$ with eigenvalue $1-(\kappa-1)!$.
		\end{enumerate}
	\end{proposition}
	
	\begin{proof}
		Assume first that $G$ is $\kappa$-uniform. By Lemma~\ref{lem:degrees}, every $\xi\in\cV$ with $v:=\out(\xi)$ satisfies
		\[
		\deg^{\mathrm{out}}_{\cG}(\xi)=(\deg v-1)\,(\kappa-1)!=(\kappa-1)!\;d(\xi).
		\]
		Hence $\cD_{\cG}=(\kappa-1)!\,\cD$, and therefore
		\[
		\cL=\Id-\cD^{-1}\cA=\Id-(\kappa-1)!\,\cD_{\cG}^{-1}\cA
		=\bigl(1-(\kappa-1)!\bigr)\Id+(\kappa-1)!\,\cL_{\mathrm{rw}} .
		\]
		
		Conversely, assume that $\cL=a\,\Id+b\,\cL_{\mathrm{rw}}$. Let $\xi\in\cV$ be a permutation of $e_i$ with output $v$. By Lemma~\ref{lem:degrees}, $\xi$ has at least one out-neighbor $\eta$ in $\cG$, and comparing the entries in position $(\xi,\eta)$ gives $1/d(\xi)=b/\deg^{\mathrm{out}}_{\cG}(\xi)$, that is,
		\[
		b=\frac{\deg^{\mathrm{out}}_{\cG}(\xi)}{d(\xi)}=\frac{\pdeg v-(c_i-1)!}{\deg v-1}.
		\]
		If two edges $e_i$ and $e_j$ share a vertex $v$, applying this to a permutation of $e_i$ and to a permutation of $e_j$, both with output $v$, gives $(c_i-1)!=(c_j-1)!$, i.e.\ $c_i=c_j$. Hence, every connected component of $G$ is uniform. Moreover, if the edges of a component have cardinality $c$, then for every vertex $v$ of this component we have $\pdeg v=\deg v\,(c-1)!$, and the identity above gives $b=(c-1)!$. Since $b$ does not depend on the component, all components have the same edge cardinality, that is, $G$ is uniform. Finally, comparing the diagonal entries, which are all equal to $1$ in both $\cL$ and $\cL_{\mathrm{rw}}$, gives $a=1-b$.
		
		Claim 1 follows from the identity $$\cL=\bigl(1-(\kappa-1)!\bigr)\Id+(\kappa-1)!\,\cL_{\mathrm{rw}},$$ which maps each eigenvalue $z$ of $\cL_{\mathrm{rw}}$ to $1+(\kappa-1)!\,(z-1)$, with the same eigenvectors. For Claim 2, observe that every node of $\cG$ has at least one out-neighbor by Lemma~\ref{lem:degrees}, so that $\cD_{\cG}^{-1}\cA$ is row-stochastic and the constant functions lie in the kernel of $\cL_{\mathrm{rw}}$. By the identity above, they are therefore eigenfunctions of $\cL$ with eigenvalue $1-(\kappa-1)!$.
	\end{proof}

	\section{The reduced non-backtracking Laplacian}\label{sec:incidences}
	
	Recall that the \emph{incidences} of $G$ are the edge--vertex pairs $(e_i,v)$ with $v\in e_i$, and let
	\[
	I(G):=\bigl\{(e_i,v)\,:\, v\in e_i\bigr\},
	\quad\text{so that}\quad
	|I(G)|=\sum_{i=1}^M c_i .
	\]
	In this section we show that, although $\cL$ is defined on the $\sum_i c_i!$ permutations of the edges, its spectrum can be computed from a much smaller matrix indexed by the $\sum_i c_i$ incidences of $G$, namely the reduced non-backtracking Laplacian of Definition~\ref{def:reduced} below. For instance, for the laminar hypergraph $\Lambda_{2,3}$ of Example~\ref{ex:laminar}, this reduces the size of the eigenvalue problem from $40376$ to $24$. This is possible because, by definition of the non-backtracking graph, the out-neighbors of a permutation in $\cG$ depend only on its edge and its output, and its in-neighbors depend only on its edge and its input. As a consequence, the eigenvalues of $\cL$ are those of the reduced non-backtracking Laplacian, together with $\sum_i(c_i!-c_i)$ additional eigenvalues that are equal to $1$ (Theorem~\ref{thm:reduction}). The permutations still play a role in the spectrum, but only through their number: the entry of the reduced matrix corresponding to a step through $e_j$ from $v$ to $w$ is proportional to $(c_j-2)!$, which is the number of permutations of $e_j$ with input $v$ and output $w$.
	
	Since the spectra of $\cL$ and of the reduced non-backtracking Laplacian determine each other, one might take the latter as the definition of the non-backtracking Laplacian of a hypergraph. We explain in Remark~\ref{rem:whyL} why we nevertheless work with $\cL$.
	
	\begin{definition}\label{def:reduced}
		The \emph{reduced non-backtracking Laplacian} of $G$ is the $|I(G)|\times|I(G)|$ matrix $$\widehat{\cL}:=\Id-\cR,$$ where $\cR$ is the matrix defined by
		\[
		\cR_{(e_i,v),(e_j,w)} :=
		\begin{cases}
		\dfrac{(c_j-2)!}{\deg v-1}
		& \text{if } j\neq i,\ v\in e_j,\ w\in e_j \text{ and } w\neq v,\\[2ex]
		0 & \text{otherwise.}
		\end{cases}
		\]
	\end{definition}
	
	In words, the incidence $(e_i,v)$ represents a walker that has just left $e_i$ through $v$, and $\cR$ moves it through a different edge $e_j\ni v$ and out of it through a vertex $w\neq v$, with weight $(c_j-2)!/(\deg v-1)$, where $(c_j-2)!$ counts the permutations of $e_j$ with input $v$ and output $w$. 
	
	If $G$ is a graph, then $(c_j-2)!=1$ and $\widehat{\cL}$ is the non-backtracking Laplacian $\mathcal{L}$, after identifying each incidence $(e_i,v)$ with the orientation of $e_i$ whose output is $v$.

	The key ingredient in the proofs of Theorem~\ref{thm:reduction} and Corollary~\ref{cor:gap} is the following factorization of $\cD^{-1}\cA$ through the incidences. In the statement, for a positive integer $n$, we let $J_n$ denote the $n\times n$ all-ones matrix.
	
	\begin{lemma}\label{lem:factorization}
		Let $\Phi$ and $\Psi$ be the $|\cV|\times|I(G)|$ matrices defined, for $\xi\in\cV$ and $(e_j,w)\in I(G)$, by
		\begin{align*}
			&\Phi_{\xi,(e_j,w)}:=\mathbbm{1}[\xi\text{ is a permutation of }e_j]\cdot
			\mathbbm{1}\bigl[\out(\xi)=w\bigr],\\
			&\Psi_{\xi,(e_j,w)}:=\mathbbm{1}[\xi\text{ is a permutation of }e_j]\cdot
			\mathbbm{1}\bigl[\inp(\xi)=w\bigr],
		\end{align*}
		and let $T$ and $\widetilde{D}$ be the $|I(G)|\times|I(G)|$ matrices
		\begin{align*}
			&T_{(e_i,v),(e_j,w)}:=\mathbbm{1}[v=w]\cdot\mathbbm{1}[i\neq j],\\
			&\widetilde{D}:=\operatorname{diag}\bigl(\deg v-1\bigr)_{(e_i,v)\in I(G)} .
		\end{align*}
		Moreover, let $H:=\widetilde{D}^{-1}T$ and $C:=\Psi^\top\Phi$. Then:
		\begin{enumerate}
			\item $\cD^{-1}\cA=\Phi\,H\,\Psi^\top$;
			\item $\cR=H\,C$;
			\item $H$ and $C$ are symmetric. Moreover, grouping the incidences by their vertex, $H$ is block diagonal with blocks
			\[
			(J_{\deg v}-\Id_{\deg v})/(\deg v-1),\quad v\in V,
			\]
			and therefore its eigenvalues are $1$ and $-1/(\deg v-1)$. Grouping the incidences by their edge, $C$ is block diagonal with blocks
			\[
			(c_i-2)!\,(J_{c_i}-\Id_{c_i}),\quad i\in\{1,\dots,M\},
			\]
			and therefore its eigenvalues are $(c_i-1)!$ and $-(c_i-2)!$.
		\end{enumerate}
		In particular, $H$, $C$ and $\cR$ are invertible.
	\end{lemma}

	\begin{proof}
		Note that each row of $\Phi$ (respectively, of $\Psi$) contains exactly one non-zero
		entry, marking the incidence given by the output (respectively, the input)
		of the corresponding permutation. The matrix $T$ encodes the admissible
		changes of edge, since the walk leaves $e_i$ and enters a different edge
		through the same vertex, and $\widetilde{D}$ normalizes by the number of such
		choices.
		\begin{enumerate}
			\item Let $\xi,\eta\in\cV$ be permutations of $e_i$ and $e_j$, respectively. A direct computation gives
			\[
			(\Phi\,T\,\Psi^\top)_{\xi\eta}
			=\mathbbm{1}\bigl[\out(\xi)=\inp(\eta)\bigr]
			\cdot\mathbbm{1}[i\neq j]=\cA_{\xi\eta},
			\]
			so that $\cA=\Phi T\Psi^\top$. Moreover, $\cD^{-1}\Phi=\Phi\widetilde{D}^{-1}$, since $\Phi_{\xi,(e_j,w)}\neq 0$ only if $w=\out(\xi)$, in which case $d(\xi)=\deg w-1$. Hence,
			\[
			\cD^{-1}\cA=\Phi\widetilde{D}^{-1}T\Psi^\top=\Phi H\Psi^\top.
			\]
			\item Observe that
			\begin{align*}
				(\Psi^\top\Phi)_{(e_i,v),(e_j,w)}
				&=\mathbbm{1}[i=j]\cdot
				\#\bigl\{\xi\text{ permutation of }e_i:\inp(\xi)=v\text{ and }
				\out(\xi)=w\bigr\}
				\\ &=\mathbbm{1}[i=j]\cdot(c_i-2)!\cdot\mathbbm{1}[v\neq w],
			\end{align*}
			since the permutations of $e_i$ with fixed distinct input $v$ and
			output $w$ are exactly $(c_i-2)!$, while no permutation has equal input and
			output. Therefore,
			\begin{align*}
				(H C)_{(e_i,v),(e_j,w)}
				&=\frac{1}{\deg v-1}\sum_{(e_l,u)\in I(G)}
				\mathbbm{1}[u=v]\,\mathbbm{1}[l\neq i]\cdot
				\mathbbm{1}[l=j]\,(c_j-2)!\,\mathbbm{1}[u\neq w]\\
				&=\frac{(c_j-2)!}{\deg v-1}\cdot
				\mathbbm{1}[j\neq i]\,\mathbbm{1}[v\in e_j]\,
				\mathbbm{1}[w\in e_j]\,\mathbbm{1}[w\neq v]\\
				&=\cR_{(e_i,v),(e_j,w)} .
			\end{align*}
			
			\item Symmetry of $H$ and $C$ is clear from their entries, and the block
			structure follows from the computations above. Since $J_n-\Id_n$ has
			eigenvalues $n-1$ and $-1$, the eigenvalues are as stated, and they are
			non-zero because $\deg v\geq\delta\geq2$ and $c_i\geq\kappa\geq2$. Hence $H$
			and $C$ are invertible, and so is $\cR=HC$.
		\end{enumerate}
	\end{proof}

	\begin{example}\label{ex:blocks}
		We illustrate Lemma~\ref{lem:factorization} on a single edge $e_i=\{v_1,v_2,v_3\}$, with the permutations ordered as in Example~\ref{ex:threeedge}. With rows and columns grouped by edge, $\Phi$ and $\Psi$ are block diagonal, with one block for each edge. The blocks of $\Phi$ and $\Psi$ corresponding to $e_i$, with rows indexed by $e_i^{(1)},\dots,e_i^{(6)}$ and columns by the incidences $(e_i,v_1),(e_i,v_2),(e_i,v_3)$, are
		\[
		\Phi_i=\begin{pmatrix}
		0&0&1\\ 0&1&0\\ 0&0&1\\ 1&0&0\\ 0&1&0\\ 1&0&0
		\end{pmatrix},
		\qquad
		\Psi_i=\begin{pmatrix}
		1&0&0\\ 1&0&0\\ 0&1&0\\ 0&1&0\\ 0&0&1\\ 0&0&1
		\end{pmatrix}.
		\]
		Here, each row of $\Phi_i$ marks the output of the corresponding permutation, and each row of $\Psi_i$ its input. In particular, two permutations of $e_i$ have the same row of $\Phi_i$ precisely when they have the same output. Moreover, the corresponding block of $C=\Psi^\top\Phi$ is
		\[
		\Psi_i^{\top}\Phi_i=\begin{pmatrix}
		0&1&1\\ 1&0&1\\ 1&1&0
		\end{pmatrix}
		=(c_i-2)!\,\bigl(J_3-\Id_3\bigr),
		\]
		whose $(v,w)$ entry counts the permutations of $e_i$ with input $v$ and output $w$, which here is exactly one for every $v\neq w$. For an edge $e_i$ with $c_i=4$, the corresponding block would be $2\,(J_4-\Id_4)$, the factor $2=(c_i-2)!$ recording the two possible orders of the two vertices of $e_i$ other than the input and the output. Through the identity $\cR=HC$ of Lemma~\ref{lem:factorization}, this is where the factor $(c_j-2)!$ in the entries of $\cR$ comes from.
	\end{example}

	We can now prove the main result of this section: by Lemma~\ref{lem:factorization}, the non-backtracking Laplacian $\cL$ is similar to the direct sum of the reduced non-backtracking Laplacian $\widehat{\cL}$ and an identity matrix. For a multiset, we write $\{\lambda\}^{(m)}$ for the element $\lambda$ repeated $m$ times.

	\begin{theorem}\label{thm:reduction}
		The non-backtracking Laplacian $\cL$ is similar to the block diagonal matrix
		\[
		\widehat{\cL}\oplus\Id_q,\quad \text{where}\quad q:=\sum_{i=1}^M(c_i!-c_i).
		\]
		In particular:
		\begin{enumerate}
			\item as multisets, $\sigma(\cL)=\sigma(\widehat{\cL})\sqcup\{1\}^{(q)};$
			\item $1$ is an eigenvalue of $\cL$ if and only if $G$ is not a graph, in which case it is semisimple, with multiplicity $q$;
			\item $\cL$ is diagonalizable if and only if $\widehat{\cL}$ is diagonalizable.
		\end{enumerate}
	\end{theorem}

	\begin{proof}
		Let $X:=\Phi$ and $Y:=H\Psi^\top$. By Lemma~\ref{lem:factorization}, $\cD^{-1}\cA=XY$ and $YX=\cR$, which is invertible. It follows that $X$ has full column rank $|I(G)|$ and $Y$ has full row rank, therefore $\ker Y$ has dimension $|\cV|-|I(G)|=q$. Moreover, $\operatorname{im}X\cap\ker Y=\{0\}$: if $X\mathbf{u}\in\ker Y$, then $\cR\mathbf{u}=YX\mathbf{u}=0$, therefore $\mathbf{u}=0$. Hence, $\C^{|\cV|}=\operatorname{im}X\oplus\ker Y$, and both subspaces are invariant under $XY$, since $XY(X\mathbf{u})=X(\cR\mathbf{u})$ and $XY\mathbf{z}=0$ for $\mathbf{z}\in\ker Y$. In a basis formed by the columns of $X$ followed by a basis of $\ker Y$, the matrix $\cD^{-1}\cA=XY$ is therefore represented by $\cR\oplus 0_q$, and $\cL$ by $\widehat{\cL}\oplus\Id_q$. We now prove the three claims.
		\begin{enumerate}
			\item Similar matrices have the same eigenvalues, with the same algebraic multiplicities, and the eigenvalues of a block diagonal matrix are those of its blocks.
			\item Since $\cR$ is invertible, $1$ is not an eigenvalue of $\widehat{\cL}=\Id-\cR$. Hence, the eigenvalue $1$ of $\cL$ only comes from the block $\Id_q$, so that its algebraic and geometric multiplicities both equal $q$, and it is semisimple. Moreover, $q>0$ if and only if $c_i\geq 3$ for some $i$, that is, if and only if $G$ is not a graph.
			\item A block diagonal matrix is diagonalizable if and only if each of its blocks is, and $\Id_q$ is diagonal.
		\end{enumerate}
	\end{proof}

	\begin{remark}[Why $\cL$ and not $\widehat{\cL}$]\label{rem:whyL}
		By Theorem~\ref{thm:reduction}, the spectra of $\cL$ and $\widehat{\cL}$ determine each other, since $1\notin\sigma(\widehat{\cL})$ and $q$ only depends on the cardinalities of the edges. Moreover, if $G$ is a graph, then $\widehat{\cL}=\cL$, after identifying each incidence $(e_i,v)$ with the orientation of $e_i$ whose output is $v$. One might therefore take $\widehat{\cL}$ as the definition of the non-backtracking Laplacian of a hypergraph. We nevertheless work with $\cL$, because it is conceptually more natural. Indeed, $\cL$ is obtained from the $(0,1)$ adjacency matrix of a directed graph, normalized by the number of edges through which the walk can continue, whereas the weights $(c_j-2)!/(\deg v-1)$ of $\widehat{\cL}$ would be hard to motivate directly: they arise from $\cL$, as the numbers of permutations of $e_j$ with fixed input and output.
	\end{remark}

	While Theorem~\ref{thm:reduction} determines the multiplicity of the eigenvalue $1$, Lemma~\ref{lem:factorization} also allows us to bound the other eigenvalues of $\cL$ away from $1$.

	\begin{corollary}\label{cor:gap}
		Every eigenvalue $\lambda\neq 1$ of $\cL$ satisfies
		\[
		\frac{(\kappa-2)!}{\Delta-1}\;\leq\;|1-\lambda|\;\leq\;(K-1)!\,.
		\]
	\end{corollary}
	
	\begin{proof}
		By Theorem~\ref{thm:reduction}, $\lambda\in\sigma(\widehat{\cL})$, that is, $\lambda=1-\mu$ for some $\mu\in\sigma(\cR)$, so that $|1-\lambda|=|\mu|$. By Lemma~\ref{lem:factorization}, $\cR=HC$, where $H$ and $C$ are symmetric and invertible. Since the spectral norm of a symmetric matrix is the largest absolute value of its eigenvalues, the eigenvalues given in Lemma~\ref{lem:factorization} yield
		\[
		\|H\|=1,\qquad \|H^{-1}\|=\Delta-1,\qquad \|C\|=(K-1)!,\qquad \|C^{-1}\|=\frac{1}{(\kappa-2)!}.
		\]
		Since every eigenvalue of a matrix is bounded in absolute value by its spectral norm, and $\mu^{-1}$ is an eigenvalue of $\cR^{-1}=C^{-1}H^{-1}$, we obtain
		\[
		|\mu|\leq\|H\|\,\|C\|=(K-1)!
		\qquad\text{and}\qquad
		|\mu|^{-1}\leq\|C^{-1}\|\,\|H^{-1}\|=\frac{\Delta-1}{(\kappa-2)!},
		\]
		which proves the claim.
	\end{proof}

	In Corollary~\ref{cor:gap}, the upper bound  is the radius of the disk of Proposition~\ref{prop:gershgorin}, while the lower bound is new. For graphs, the lower bound is $1/(\Delta-1)$, which is the bound of~\cite[Section 4]{NB-Laplacian}.

	\begin{remark}
		Operators indexed by the incidences of a hypergraph also appear in~\cite{CEH}, where the incidences are called \emph{pointed edges}. There, the non-backtracking operator, which we denote by $\mathcal{K}$, has entry $1$ from $(v,e_i)$ to $(w,e_j)$ if and only if $w\in e_i\setminus\{v\}$, $w\in e_j$ and $e_j\neq e_i$ \cite[Definitions 2.1--2.2]{CEH}: the walk crosses $e_i$ from $v$ to $w$, and then moves to another edge at $w$. Identifying $(v,e_i)$ with the incidence $(e_i,v)$, Definition~\ref{def:reduced} gives
		\[
		\cR=\widetilde{D}^{-1}\,\mathcal{K}^{\top}\,W,
		\qquad
		W:=\operatorname{diag}\bigl((c_j-2)!\bigr)_{(e_j,w)\in I(G)},
		\]
		with $\widetilde{D}$ as in Lemma~\ref{lem:factorization}. That is, $\cR$ is obtained from $\mathcal{K}^{\top}$ by weighting each crossing of an edge $e_j$ by the number $(c_j-2)!$ of permutations of $e_j$ with the corresponding input and output, and by normalizing by the number of edges through which the walk can continue. Since $\widetilde{D}$ and $W$ are both diagonal, conjugating by $W$ shows that $\cR$ is similar to $W\widetilde{D}^{-1}\mathcal{K}^{\top}$, where
		\[
		W\widetilde{D}^{-1}
		=\operatorname{diag}\left(\frac{(c_i-2)!}{\deg v-1}\right)_{(e_i,v)\in I(G)} .
		\]
		In particular, if $(c_i-2)!/(\deg v-1)$ takes the same value $\nu$ on all incidences, as is the case when all vertices have the same degree $\Delta$ and either $G$ is $\kappa$-uniform, with $\nu=(\kappa-2)!/(\Delta-1)$, or all edges have at most $3$ vertices, with $\nu=1/(\Delta-1)$, then $\sigma(\cR)=\nu\,\sigma(\mathcal{K})$. In general, $W\widetilde{D}^{-1}$ is not scalar, and the two spectra need not be related in such a simple way.
	\end{remark}

	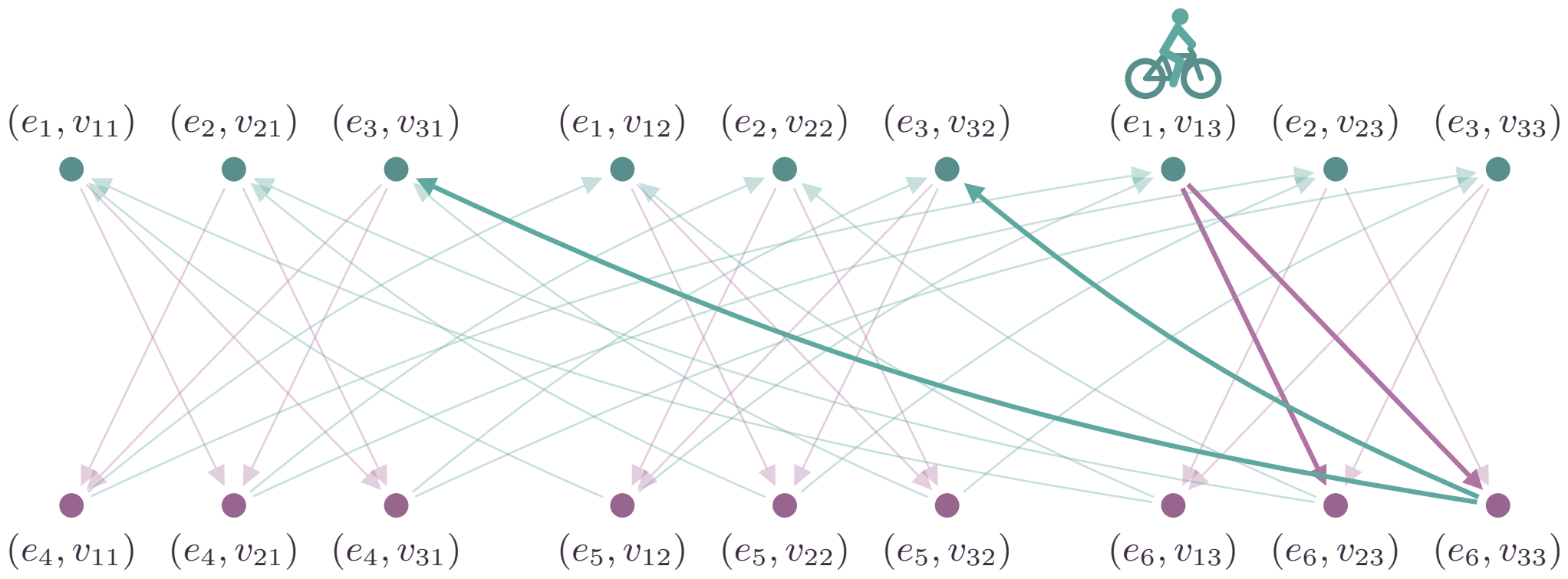
\begin{figure}[h]\centering
		\tikzset{
			incr/.style={circle, fill=color1!80!nodec, inner sep=0pt, minimum size=6pt},
			incc/.style={circle, fill=color2!80!nodec, inner sep=0pt, minimum size=6pt},
			arcc/.style={-{Triangle[length=5pt,width=4.5pt]}, line width=1.2pt, color2,
				shorten >=2pt, shorten <=2pt},
			arcr/.style={-{Triangle[length=5pt,width=4.5pt]}, line width=1.2pt, color1,
				shorten >=2pt, shorten <=2pt},
			faint/.style={line width=0.5pt, opacity=0.35},
		}
		\begin{tikzpicture}
		\foreach \a in {1,2,3} \foreach \b in {1,2,3}{
			\pgfmathsetmacro{\x}{(\b-1)*4.75+(\a-1)*1.4}
			\node[incr] (r\a\b) at (\x,2.9) {};
			\node[incc] (c\a\b) at (\x,0) {};
			\node[lab,nodec,font=\footnotesize,above=3pt] at (r\a\b) {$(e_{\a},v_{\a\b})$};
			\pgfmathtruncatemacro{\j}{3+\b}
			\node[lab,nodec,font=\footnotesize,below=3pt] at (c\a\b) {$(e_{\j},v_{\a\b})$};}
		\foreach \a in {1,2,3} \foreach \b in {1,2,3} \foreach \aa in {1,2,3}{
			\ifnum\a=\aa\else
			\ifnum\a\b=13
			\draw[arcc] (r\a\b) -- (c\aa\b);
			\else
			\draw[arcc,faint] (r\a\b) -- (c\aa\b);
			\fi\fi}
		\foreach \a in {1,2,3} \foreach \b in {1,2,3} \foreach \bb in {1,2,3}{
			\ifnum\b=\bb\else
			\ifnum\a\b=33
			\draw[arcr] (c\a\b) to[bend left=8] (r\a\bb);
			\else
			\draw[arcr,faint] (c\a\b) to[bend left=8] (r\a\bb);
			\fi\fi}
		\cyclist[0.6]{($(r13)+(0,0.62)$)}{0}
		\end{tikzpicture}
		\caption{The directed graph on the $18$ incidences of the $3$-lattice whose weighted
			adjacency matrix is $\cR$; all weights equal $(3-2)!/(2-1)=1$. Top:
			the incidences of the rows; bottom: those of the columns. From $(e_a,v_{ab})$,
			the walk crosses the column $e_{3+b}$ and exits it at either of its two other
			vertices (in purple); from $(e_{3+b},v_{ab})$, it crosses the row $e_a$ and exits
			it at either of its two other vertices (in teal). The arcs out of
			$(e_1,v_{13})$ and $(e_6,v_{33})$ are highlighted: they continue the walk of
			Figure~\ref{fig:lattice}, from the cyclist's position $(e_1,v_{13})$.}
		\label{fig:lattice-incidences}
	\end{figure}

	\begin{example}[Example~\ref{ex:lattice}, continued]\label{ex:lattice2}
		Let $G$ be the $3$-lattice of Example~\ref{ex:lattice}. Its set of incidences $I(G)$ has $\sum_i c_i=18$ elements and, since $\deg v=2$ and $c_j=3$ for all $v$ and $j$, all non-zero entries of $\cR$ equal $(c_j-2)!/(\deg v-1)=1$: that is, $\cR$ is the adjacency matrix of the directed graph on $I(G)$ in which, from an incidence $(e_i,v)$, the walk moves to the unique other edge $e_j\ni v$ and exits it at either of the two vertices $w\neq v$ of $e_j$ (Figure~\ref{fig:lattice-incidences}). A direct computation shows that the characteristic polynomial of $\cR$ is
		\[
		\bigl(x^2-4\bigr)\bigl(x^2-1\bigr)^4\bigl(x^2+2\bigr)^4 ,
		\]
		so that $\sigma(\cR)$ consists of $\pm 2$, of $\pm 1$ with multiplicity $4$ each, and of $\pm i\sqrt{2}$ with multiplicity $4$ each. Hence, by Theorem~\ref{thm:reduction}, as a multiset,
		\[
		\sigma(\cL)
		= \{1\}^{(18)}
		\sqcup \{-1,\,3\}
		\sqcup \{0\}^{(4)}
		\sqcup \{2\}^{(4)}
		\sqcup \bigl\{1-i\sqrt{2}\bigr\}^{(4)}
		\sqcup \bigl\{1+i\sqrt{2}\bigr\}^{(4)} .
		\]
		We point out three facts. First, $-1=1-(\kappa-1)!$ is the eigenvalue of the constant functions, in accordance with Proposition~\ref{prop:uniform}. Second, $\sigma(\cL)$ is symmetric with respect to the map $\lambda\mapsto 2-\lambda$. As in the graph case~\cite[Proposition 3.4]{NB-Laplacian}, this follows from the bipartiteness of $\cG$ observed in Example~\ref{ex:lattice}. Third, both bounds of Corollary~\ref{cor:gap} are sharp in this example: the upper bound $(K-1)!=2$, which is also the radius of the disk of Proposition~\ref{prop:gershgorin}, is attained by the eigenvalues $-1$ and $3$, and the lower bound $(\kappa-2)!/(\Delta-1)=1$ by the eigenvalues $0$ and $2$.
	\end{example}

	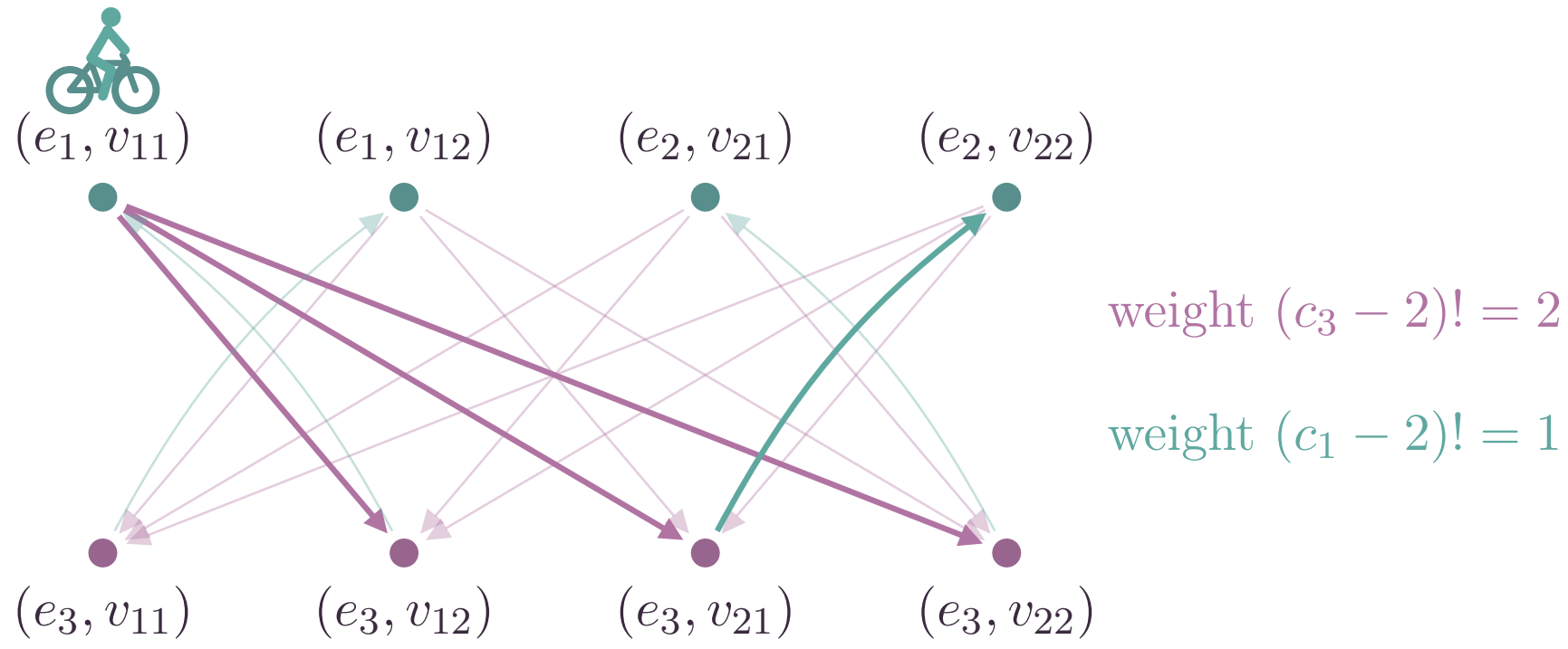
\begin{figure}[htb]\centering
		\tikzset{
			incr/.style={circle, fill=color1!80!nodec, inner sep=0pt, minimum size=6pt},
			incc/.style={circle, fill=color2!80!nodec, inner sep=0pt, minimum size=6pt},
			arcc/.style={-{Triangle[length=5pt,width=4.5pt]}, line width=1.2pt, color2,
				shorten >=2pt, shorten <=2pt},
			arcr/.style={-{Triangle[length=5pt,width=4.5pt]}, line width=1.2pt, color1,
				shorten >=2pt, shorten <=2pt},
			faint/.style={line width=0.5pt, opacity=0.35},
		}
		\begin{tikzpicture}[scale=1]
		\foreach \v/\x in {11/0,12/2.2,21/4.4,22/6.6}{
			\node[incr] (t\v) at (\x,2.6) {};
			\node[incc] (b\v) at (\x,0) {};}
		\node[lab,nodec,above=3pt] at (t11) {$(e_1,v_{11})$};
		\node[lab,nodec,above=3pt] at (t12) {$(e_1,v_{12})$};
		\node[lab,nodec,above=3pt] at (t21) {$(e_2,v_{21})$};
		\node[lab,nodec,above=3pt] at (t22) {$(e_2,v_{22})$};
		\foreach \v in {11,12,21,22}{\node[lab,nodec,below=3pt] at (b\v) {$(e_3,v_{\v})$};}
		\foreach \v in {11,12,21,22}{
			\foreach \w in {11,12,21,22}{
				\ifnum\v=\w\else
				\ifnum\v=11
				\draw[arcc] (t\v) -- (b\w);
				\else
				\draw[arcc,faint] (t\v) -- (b\w);
				\fi
				\fi}}
		\draw[arcr,faint] (b11) to[bend left=12] (t12);
		\draw[arcr,faint] (b12) to[bend right=12] (t11);
		\draw[arcr]       (b21) to[bend left=12] (t22);
		\draw[arcr,faint] (b22) to[bend right=12] (t21);
		\node[lab,color2,anchor=west] at (7.2,1.75) {weight $(c_3-2)!=2$};
		\node[lab,color1,anchor=west] at (7.2,0.85) {weight $(c_1-2)!=1$};
		\cyclist[0.6]{($(t11)+(0,0.62)$)}{0}
		\end{tikzpicture}
		\caption{The weighted directed graph on the $8$ incidences of $\Lambda_{2,2}$
			whose weighted adjacency matrix is $\cR$. From an incidence of a level-$1$
			edge, the walk crosses $e_3$ and exits it at any of its three other vertices
			(in purple, weight $2$); from an incidence $(e_3,v)$, it crosses the level-$1$
			edge of $v$ and exits it at its other vertex (in teal, weight $1$). The arcs out of $(e_1,v_{11})$ and $(e_3,v_{21})$ are highlighted: they show two steps of the walk from the cyclist's position $(e_1,v_{11})$.}
		\label{fig:laminar-incidences}
	\end{figure}

	\begin{example}[Example~\ref{ex:laminar}, continued]\label{ex:city2}
		Let $G=\Lambda_{2,2}$, with vertices $v_{11},v_{12},v_{21},v_{22}$, level-$1$ edges $e_1=\{v_{11},v_{12}\}$ and $e_2=\{v_{21},v_{22}\}$, and level-$2$ edge $e_3=V$. As for the $3$-lattice, $\deg v=2$ for every vertex, therefore $\cD=\Id$ and $\cL=\Id-\cA$, but in this case, $G$ is not uniform. The matrix $\cR$ has size $8$, while $\cL$ has size $|\cV|=28$ (Figure~\ref{fig:laminar-incidences}): from an incidence of a level-$1$ edge, the walk enters $e_3$ and exits it at any of its three other vertices, with weight $(c_3-2)!=2$; from an incidence $(e_3,v)$, it enters the level-$1$ edge of $v$ and exits it at its other vertex, with weight $1$. A direct computation gives the characteristic polynomial $\bigl(x^2-6\bigr)\bigl(x^2+2\bigr)\bigl(x^2-2\bigr)^2$ of $\cR$, hence, by Theorem~\ref{thm:reduction},
		\[
		\sigma(\cL)
		= \{1\}^{(20)}
		\sqcup \bigl\{1\pm\sqrt{6}\bigr\}
		\sqcup \bigl\{1\pm i\sqrt{2}\bigr\}
		\sqcup \bigl\{1\pm\sqrt{2}\bigr\}^{(2)} .
		\]
		As for the lattice, $\sigma(\cL)$ is symmetric with respect to $\lambda\mapsto 2-\lambda$, since every walk on $\cG$ alternates between the level-$1$ edges and $e_3$, therefore $\cG$ is bipartite. Unlike for the lattice, however, neither bound of Corollary~\ref{cor:gap} is sharp: the eigenvalues of $\cL$ different from $1$ lie at distance between $\sqrt{2}$ and $\sqrt{6}$ from $1$, while the two bounds equal $1$ and $(K-1)!=6$. Moreover, the constant functions are not eigenfunctions of $\cL$, since the row sums of $\cA$ equal $6$ on the permutations of the level-$1$ edges and $1$ on those of $e_3$ (Remark~\ref{rem:notprob}); accordingly, unlike in the uniform case of Proposition~\ref{prop:uniform}, $\cL$ has no eigenvalue of the form $1-(\kappa-1)!$.
	\end{example}

	\section*{Tool and computational resource disclosure}
	The authors used Claude as a tool for language polishing, for identifying possible inconsistencies or errors in preliminary drafts, and for drafting figures.
	
	\section*{Funding}
	Raffaella Mulas is supported by the Dutch Research Council (NWO) through the grant VI.Veni.232.002 (Grant ID: \url{https://doi.org/10.61686/NLUFV33835}), which also supported Eleonora Andreotti's research visit to VU Amsterdam in June 2026, during which this project was developed. Eleonora cycled from Trento to Amsterdam for this visit and returned by bus,
	so her bike ride was non-backtracking.
	
	\bibliographystyle{plain} 
	
	\bibliography{Bibliography}
\end{document}